\documentstyle [twoside]{article} 
\newtheorem{thm}{Th\'eor\`eme}[section]
\newtheorem{prop}[thm]{Proposition}
\newtheorem{cor}[thm]{Corollaire}
\newtheorem{lem}[thm]{Lemme}
\newtheorem{remark}{Remarque} 
\newtheorem{dfn}{D\'efinition} 
\def\Blat{\mbox{\it \raise2pt\hbox{"}\kern-2pt H}}
\def\lvup{\rlap{\ ${}^{\ell\atop{\hbox{${}^{\vee}$}}}$}\cdots} 
\begin{document} 
\begin{center} 
{\center{\Large{\bf 
Transform\'ee  de Mellin des int\'egrales- fibres \\ 
associ\'ees  aux   singularit\'es isol\'ees \\
 d'intersection compl\`ete quasihomog\`enes} }} 
 
 \vspace{1pc} 
{ \center{\large{ Susumu TANAB\'E }}} 
\end{center} 
 
\noindent 
\begin{center} 
 \begin{minipage}[t]{10.2cm} 
{\sc R\'esum\'e.} {\em La transform\'ee de Mellin de  l'int\'egrale 
-fibre est calcul\'ee pour certaines singularit\'es isol\'ees d'intersection
compl\`ete quasihomog\`enes (surtout singularit\'es unimodales
de la liste de Giusti et de Wall). On montre 
la propri\'et\'e de la sym\'etrie des spectres de Gauss-Manin 
(Th\'eor\`eme ~\ref{thm1}) et on met \`a jour
la structure de r\'eseaux des p\^oles de 
la transform\'ee de Mellin, exprim\'ee au moyen des donn\'ees 
topologiques des singularit\'es (Th\'eor\`eme ~\ref{thm43}, Th\'eor\`eme 
~\ref{thm52}). Comme application de ces r\'esultats, on exprime le nombre 
de Hodge de la fibre par le d\'enombrement de  spectres de Gauss-Manin.} 
 \end{minipage} \hfill 
\end{center} 
 \vspace{1pc} 
{ 
\center{\section{Introduction}} 
} 
 
Ici on calcule concr\`etement  le  syst\`eme  de  Gauss-Manin  de 
l'int\'egrale-fibre associ\'ee aux singularit\'es isol\'ees 
 d'intersection compl\`ete (SIIC) quasihomog\`ene. 
 
D'abord nous fixons la situation. 
Pour les deux vari\'et\'es complexes $ X = ( {\bf C}^{n+k}, 0), 
S = ( {\bf C}^k, 0), $ on regarde une application quasihomog\`ene 
d'intersection compl\`ete, 
$$ f: X \rightarrow S  \leqno(0.1)$$ 
telle que 
$X_s :=\{(x_1,\cdots, x_{n+k})\in X ; 
f_1(x) = s_1, \ldots, f_k(x) = s_k 
\}.$ 
C'est \`a dire que $f_1(x), \cdots,f_k(x)$ sont 
des polyn\^omes quasihomog\`enes par rapport \`a un poids
et $dim X_0 =n \geq 1.$ On suppose en  plus 
que   $f$   poss\`ede   une   singularit\'e    isol\'ee \`a    l'origine 
i. e. $df(x)=0$ 
pour $x \in X_0$ si et seulement si $x = 0.$ 
 
Dans ce travail, notre but est de d\'ecrire les solutions explicites
du syst\`eme de
Gauss-Manin associ\'e \`a SIIC quasihomog\`enes 
pour certains  cas de courbe espace i.e. $n=1, k =2.$
Ici nous nous servons de la transform\'ee de Mellin 
d'int\'egrales-fibres parce que elle permet de mieux visualiser
les propri\'et\'es importantes de singularit\'es
SIIC.
Cette situation a
incit\'e  certains auteurs comme C.Sabbah \cite{Sab2}, \cite{Sab3},
D.Barlet\cite{Bar}, 
F.Loeser \cite{LS} \`a poursuivre
des recherches sur la  transform\'ee de Mellin d'int\'egrales fibres 
qui etaient, entre autres, motiv\'es par une id\'ee de P.Deligne
reproduite dans \cite{Sab1}.

Le plan de cet article est le suivant.
Dans le \S 1, selon Greuel et Hamm \cite{GH}, 
on introduit les espaces vectoriels  sur lesquels le syst\`eme de 
Gauss-Manin associ\'e \`a SIIC quasihomog\`ene de dimension 
arbitraire sera d\'efini et repr\'esent\'e au moyen des matrices implicitement
d\'efinies. 
Les r\'esultats du \S 2  montrent le calcul concret du syst\`eme de 
Gauss-Manin associ\'e aux 
 singularit\'es  isol\'ees simples 
d'intersection compl\`ete (SISIC) $S_{\mu},$ $U_{\mu}, T_{\mu}, W_{\mu}, 
Z_{\mu}\;\;(\mu \geq 5)$ de la liste de M.Giusti 
\cite{Gui}. Je tiens \`a noter que le calcul effectu\'e par S.Guzev \cite{Guz}
sert \`a l'\'etablissement de r\'esultats de cette section.  
%
 \vspace{1.5pc}
\footnoterule

\footnotesize{AMS Subject Classification: 14M10, 32S25, 32S40.

Key words and phrases: Gauss-Manin connexion, complete intersection,
Hodge structure.}

 \footnotesize{${}^1$  Travail   r\'ealis\'e   par   le   soutien
financier d'homme d'affaires M.Mikhail S.Gavounas (Moscou, Russie) et du 
Max Planck
Institut f\"ur Mathematik.
}
\normalsize
\newpage
\par

Dans le \S 3, les spectres du   syst\`eme de Gauss-Manin associ\'e \`a SISIC 
 sont d\'efinis et 
la sym\'etrie entre eux est \'etablie.  On note que 
la sym\'etrie des spectres de la  structure  de  Hodge 
mixte sur la cohomologie relative d'une SIIC 
a \'et\'e d\'emontr\'ee par W.Ebeling et J.Steenbrink 
\cite{ESt}. Notre approche est diff\'erent de celui de Ebeling-Steenbrink,
puisque nos objets principaux sur lesquel la transformation de monodromie
agit sont les espaces $V$ et $\Phi$ de Greuel-Hamm. 
 Dans le \S 4, la transform\'ee de Mellin de l'int\'egrale- fibre est 
d\'ecrite 
 au moyen des spectres mentionn\'es. Pour cela, on r\'esout une \'equation
aux diff\'erences finies. L'interpr\'etation de 
l'int\'egrale- fibre comme fonction hyperg\'eom\'etrique g\'en\'eralis\'ee 
au sens de Mellin-Barnes-Pincherle est
donn\'ee. Les p\^oles de la transform\'ee de Mellin donnent des informations 
sur la $b-$ fonction en 2-variables introduite par Sabbah \cite{Sab2}. 
 Dans le \S 5, en se servant du  
caract\`ere  assez universel des calculs pour SISIC,
on   g\'en\'eralise des r\'esultats des \S 2 et \S 3 
aux s\'eries des singularit\'es  "non-resonantes" et unimodales.
 Dans le \S 6, on exprime le nombre de Hodge de la fibre de Milnor 
$h^{pq}(X_s)$
au moyen des spectres de Gauss-Manin obtenus dans \S3, \S5.

Je tiens \`a remercier D.Barlet, E.Brieskorn, J.H.M.Steenbrink et  C.Sabbah de
 leurs critiques  
utiles et V. P. Palamodov d'avoir mis \`a ma disposition 
une copie de manuscrit \cite{Guz}.

 \vspace{2pc} 
{ 
\center{\section{Les espaces vectoriels  $V$ et $\Phi$ de Greuel-Hamm}} 
} 

{\bf 1.1} On reprend la situation et les notations de \S 0.
Dans cette section, on pr\'epare quelques lemmes sur 
l'intersection compl\`ete 
$X_s :=\{(x_1,\cdots, x_{m})\in X; f_1(x) = s_1, \ldots, f_k(x) = s_k 
\},$ de dimension $n = m-k \geq 0$
qui est d\'efinie par une collection de polyn\^omes quasihomog\`enes 
$f_1(x), \ldots, f_k(x).$ 

D'abord on commence par munir  nos objets  des poids  quasihomog\`enes. 
D\`es que $f_1,\cdots, f_k$ sont des polyn\^omes quasihomog\`enes, on peut 
attribuer 
aux 
variables $x_1,\cdots, x_{m}$ les poids  quasihomog\`enes. 
Notons les poids  de ces 
variables 
par 
$$w(x_1)=w_1, \cdots, w(x_{m})= w_{m},\leqno(1.1.1)$$ 
o\`u $w_1, \cdots, w_{m}$ sont les entiers 
positifs de 
pgdc \'egal \`a 1. 
On utilisera la notation $w(f_1)=p_1 
,$ $\cdots,$ $w(f_k)= p_k,$ en sorte que $p_1 \geq p_2 \geq \cdots 
\geq p_k.$ Il est naturel de d\'efinir le champ d'Euler 
$$ E= w_1x_1 \frac{\partial}{\partial x_1}+  \cdots 
+ w_m x_m \frac{\partial}{\partial x_m},\leqno(1.1.2)$$ 
de telle sorte que 
$$ E(f_j) = p_j f_j \;\;\;\mbox{pour}\;j=1,\cdots, k.$$ 
On peut associer \`a une fonction ou une forme holomorphe
 quasihomog\`ene  $\xi$ son poid quasihomog\`ene et on le note par $w(\xi).$

{\bf 1.2}
Pour calculer le syst\`eme de Gauss-Manin associ\'e aux singularit\'es 
not\'ees ci-dessus, nous introduisons les deux espaces vectoriels 
$V$ et $F,$ 
$$ F : = 
\frac{\Omega^{n+1}_X}{df_1\wedge\Omega^{n}_X +\cdots +  
df_k\wedge\Omega^{n}_X + 
i_{E}(\Omega^{n+2}_X )},\leqno(1.2.1)$$ 
o\`u $i_{E}$  signifie contraction avec le champ d'Euler $E$ d\'efini 
par (1.1.2). 
$$ V : = \frac{\Omega^{n}_X} { df_1\wedge\Omega^{n-1}_X +\cdots +  
df_k\wedge\Omega^{n-1}_X +  d\Omega^{n-1}_X  + f_1 \Omega^{n}_X + \cdots + 
f_k \Omega^{n}_X  }. \leqno(1.2.2)$$ 
L'espace      $V$ a \'et\'e, par exemple, introduit     par Greuel-Hamm 
\cite{GH}. Ils 
s'en servirent afin de calculer le nombre de Milnor $\mu$ et le 
polyn\^ome caract\'eristique de la monodromie de Picard-Lefschetz pour $f,$ 
une singularit\'e isol\'ee d'intersection compl\`ete quasihomog\`ene. 
Du lemme 3.6 de \cite{GH} on d\'eduit  que $rang_{{\bf C}}  V$  est 
\'egal  au nombre  de  Milnor  $\mu$   de   la 
singularit\'e.  Dans leur formule, $ P_V(1)= \mu$ pour le cas $n
>0.$ 
Le Satz 3.1 de \cite{GH} donne  la s\'erie  de 
Poincar\'e $P_V(t),$  
$$ P_V(t)= Res_{\tau =0} \frac{\tau^{-m+k-1}}{\tau+1}[\prod_{i=1}^m
\frac{1+ \tau t^{w_i}}{1- t^{w_i}}\prod_{j=1}^k
\frac{1 - t^{p_j}}{1+ \tau t^{p_j}} + \tau]. \leqno(1.2.3)$$
Quant \`a l'espace $F,$ on doit sa d\'efinition essentiellement \`a
S.Guzev \cite{Guz}.

{\bf 1.3}
Par les propositions suivantes, on voit l'utilit\'e de l'espace  $F$ 
pour le calcul de Gauss-Manin. Introduisons un autre espace vectoriel 
$\Phi:$ 
$$ \Phi: = 
\frac{\Omega^{m}_X}{df_1\wedge \cdots \wedge df_k \wedge \Omega^{n}_X + 
f_1\Omega^{m}_X + \cdots + f_k\Omega^{m}_X }.$$ 
On introduit un espace vectoriel $\tilde \Phi$ qui est \'evidemment
isomorphe \`a $\Phi:$ 
$$ \tilde \Phi = {\cal O}_X/\langle{\cal J}(f),f_1, \cdots, f_k \rangle,$$
o\`u ${\cal J}(f)$ l'ideal jacobien des mineurs d'ordre $k:$
$${\cal J}(f):= \langle \frac{\partial (f_1, \cdots, f_k)}
{\partial (x_{i_1},\cdots, x_{i_k})}, 1\leq i_1 < \cdots <i_k \leq n+k\rangle. $$
Pour $\phi_j(x) \in \tilde \Phi,$ on a $\phi_j (x)dx \in \Phi.$
Nous notons par $Cr(f):= \{x\in X; df_1(x) \wedge 
\cdots \wedge df_k(x)=0 \},$  l'ensemble d\'efini par l'id\'eal 
${\cal J}(f).$
 
Selon la construction de Brieskorn-Greuel \cite{Gr}, introduisons un module
$$\Blat 
= f_{\ast} \Omega^{m}_X / df_1 \wedge \cdots \wedge df_k \cdot  
d(f_{\ast}\Omega^{n-1}_X) \cong 
\frac{\Omega^{m}_X}{df_1\wedge \cdots \wedge df_k \wedge d\Omega^{n-1}_X},$$ 
qui est identifi\'e 
\`a un $ {\cal O}_S-$ module du rang $\mu$ (Proposition 2.6 \cite{Gr}). 
En fait, 
\begin{lem} 
Si $f_1,\cdots, f_k$ sont des polyn\^omes quasihomog\`enes
qui d\'efinissent SIIC,
l'espace vectoriel 
$\Phi$ 
est isomorphe \`a un autre espace vectoriel $ \Blat/(f_1,\cdots,f_k),$
le r\'eseau de Brieskorn. 
\label {lem0}
\end{lem}
{\bf D\'emonstration}
Prenons  $\alpha$ un \'el\'ement non nul  de 
$ \Omega_X^m$ admettant la d\'ecomposition
$$\alpha = \omega + \sum^k_{i=1} f_i \varphi_i +  
df_1\wedge \cdots \wedge df_k 
\wedge d\psi, $$ avec $\psi \in \Omega^{n-1}_X, \varphi_i \in \Omega^{m}_X,$
pour
$$ \omega =   df_1\wedge \cdots \wedge df_k 
\wedge \phi, $$ avec $\phi \in \Omega^{n}_X.$

Alors, 
$$\alpha = \sum^k_{i=1} f_i \varphi_i +  
df_1\wedge \cdots \wedge df_k \wedge (d\psi+ \phi). $$
L'\'enonc\'e du lemme se r\'eduit \`a la nullit\'e de $\alpha $
dans $ \Blat/(f_1,\cdots, f_k).$On peut supposer la d\'ecomposition de la forme $\phi\in \Omega^{n}_X $
selon le poids quasihomog\`ene
$$ \phi = \sum^L_{j=1}w(\phi_j)^{-1}(di_E + i_Ed)(\phi_j)= \sum^L_{j=1}\phi_j$$
o\`u $w(\phi_j) = \frac{c_1j +c_0}{c}$ avec $c, c_0, c_1,$ les 
entiers strictement positifs. Avec cette notation, la forme $\alpha$
s'\'ecrit 
$$\alpha = \sum^k_{i=1} f_i \varphi_i +  
df_1\wedge \cdots \wedge df_k \wedge [(d\psi + \sum^L_{j=1}w(\phi_j)^{-1}
(di_E + i_Ed)\phi_j]$$
$$ \equiv df_1\wedge \cdots \wedge df_k \wedge [
d( \psi +\sum^L_{j=1}w(\phi_j)^{-1}i_E(\phi_j) ) + i_E(\sum^L_{j=1}
w(\phi_j)^{-1}
d\phi_j)]  $$
$$ \equiv \sum^L_{j=1}w(\phi_j)^{-1} df_1\wedge \cdots \wedge df_k 
\wedge   
i_E(d\phi_j),$$
$$ \equiv (-1)^{k+1} i_E[\sum^L_{j=1}w(\phi_j)^{-1} df_1\wedge \cdots \wedge 
df_k \wedge   d\phi_j]\equiv 0\;\;\; \mbox{dans}\;\;  \Blat/(f_1,\cdots,f_k),$$
car $df_1\wedge \cdots \wedge df_k \wedge   d\phi_j \in \Omega^{m+1}_X
\cong 0.$ {\bf C.Q.F.D.}

De ce lemme  il suit que $rang_{\bf C}  \Phi  =  \mu.$ 
Nous notons 
les \'el\'ements de la base de $\Phi$ par $\phi_j(x)dx, 1 \leq j \leq 
\mu.$ Ici et par la suite on utilise la notation $x = (x_1, \cdots, x_m), dx
=  dx_1 \wedge \cdots \wedge dx_m. $
Notons aussi la base de l'espace $F$ par $\tilde{\omega_i}, 1 \leq i \leq 
\mu.$ 
 
  Nous soulignons ici le caract\`ere topologiquement invariant 
des espaces $F$ et $\Phi.$

${\bf 1.4.}$ 
\begin{prop}  
Si on d\'efinit le champ d'Euler comme (1.1.2), 
l'application $i_E$   don\'ee par la contraction avec $E,$ 
$$ i_E : F \rightarrow V $$ 
induit  un isomorphisme entre les deux espaces vectoriels $F$ et $V$. 
\label{prop2} 
\end{prop} 
  
{\bf D\'emonstration} 
 
1) Surjectivit\'e.  Si on prend une forme quasihomog\`ene 
$\omega \in V$, 
en vertu de la quasihomog\'en\'eit\'e de $\omega$ , 
$$ w(\omega) \omega  = i_E (d\omega ) + d( i_E \omega).    $$ 
Cela  veut dire,  pour $ \frac {d\omega}{w(\omega)} \in F ,  $ 
$$ i_E (\frac {d\omega}{w(\omega)} ) \equiv \omega \;\;\mbox {dans} \;V.$$ 
Ici on note $w(\omega)$ le poids  de la forme $\omega.$ 
 
2) Injectivit\'e. 
Supposons pour $\omega \in F,$ 
$$ i_E (\omega) = df_1 \wedge \phi_1 + \cdots + df_k \wedge \phi_k +  d\psi 
+ f_1 \omega_1 +  \cdots + f_k \omega_k, $$ 
avec $\phi_1,\cdots,   \phi_k ,  \psi \in {\Omega }^{n-1}_X,$ $  \omega_1 , 
\cdots, \omega_k \in \Omega_X^n. $ 
\c{C}a veut dire,  $$ d i_E ({\omega})= df_1 
\wedge (d\phi_1 + \omega_1) +  \cdots + df_k  \wedge 
(d\phi_k +  \omega_k) + f_1 d\omega_1 + \cdots + f_k d\omega_k$$ 
$$={w(\omega)} \omega - 
i_E (d \omega ). $$ 
Ou bien 
$$  \omega  \equiv   \frac  {1}{w(\omega)}  (f_1  d\omega_1 + \cdots +  f_k 
d\omega_k)\;\; \mbox {dans} \;\; F.  $$ 
D'autre part, puisque 
$$ p_1 f_1 d\omega_1= i_E ( df_1  \wedge  d\omega_1) + df_1 \wedge 
i_E ( d\omega_1 ), $$ 
on a 
$$ f_1 d\omega_1 \equiv 0\;\;\; \mbox {dans} \;\; F.  $$ 
D'une fa\c{c}on analogue, 
     $$ f_i d\omega_i \equiv 0 \;\;\mbox {dans} \;\; F, 2 \leq i \leq k.  $$ 

{\bf C.Q.F.D.} 

{\bf 1.5} 
En tenant compte de la  Proposition ~\ref{prop2}, 
nous notons une base de $V$ par $ \omega_i$ 
telle que $ \omega_i= i_E (\tilde{\omega}_i), \;\;1 \leq i \leq 
\mu.$ 
Dans la suite on entend par  $\tilde{\omega_i}$
une forme concr\`ete quasihomog\`ene telle que 
$(di_E +i_E d)(\tilde{\omega_i}) = \ell_i \tilde{\omega_i},$
avec le poids quasihomog\`ene   $\ell_i=w({\omega_i})$
qui figure dans les termes de la s\'erie de Poincar\'e (1.2.3), $P_F(t) = 
P_V(t).$ 
On se sert de la m\^eme convention pour la base $\phi_j(x)dx\in 
\Phi.$ 
C'est \`a dire $\tilde{\omega_i}$ est une forme 
repr\'esentant une classe d'\'equivalence, pas une classe d'\'equivalence
elle m\^eme. 

\begin{prop} 
 Pour chaque $\tilde{\omega_i},$  on a la d\'ecomposition suivante: 
$$ \tilde{\omega_i}\wedge  df_1 \wedge \lvup \wedge df_k 
\equiv \sum_{j=1}^{\mu} P^{(\ell)}_{ij}(f) 
\phi_j dx \; mod (df_1 \wedge \cdots 
\wedge df_k \wedge d\Omega^{n-1}_X ),\;\;\;\leqno(1.5.1) $$ 
avec $P^{(\ell)}_{ij}(f)\in {\bf C}[f_1,\cdots, f_k]$ et $\phi_j(x)dx\in 
\Blat/(f_1, \cdots, f_k) \cong \Phi,$
pour $1 \leq i,j \leq \mu, 1 \leq \ell \leq k$ et  
$df_1 \wedge \lvup \wedge df_k = 
\bigwedge^k_{i \not = \ell} df_i.$
\label{prop1} 
\end{prop} 
 
{\bf D\'emonstration} 
 D'apr\`es la condition d'intersection  compl\`ete sur $f$,
pour chaque $\alpha \in \Omega^{n+k},$
il existe la d\'ecomposition:
$$ 
\alpha= P(f_1, \cdots ,f_k)dx + df_1 \wedge \cdots 
\wedge df_k \wedge \beta,$$
pour certain polyn\^ome  $P(s_1, \cdots ,s_k) \in {\bf C}[s_1, \cdots s_k]$
et $\beta \in  d\Omega^{n-1}_X.$
L'unicit\'e de la d\'ecomposition d\'ecoule du fait que 
$\Blat$ est un ${\cal O}_S-$ module libre de 
rang $\mu$  engendr\'e des g\'en\'erateurs finis
(Korollar 4.9, \cite{Gr}). On applique ce raisonnement \`a la  forme
$\alpha =\tilde{\omega_i}\wedge  df_1 \wedge \lvup \wedge df_k. $
La conclusion se d\'eduit imm\'ediatement de l'isomorphisme entre 
$\Blat/(f_1, \cdots, f_k)$ et $\Phi,$
{\bf C.Q.F.D.}
 
${\bf 1.6} $
Nous abordons le calcul du syst\`eme de Gauss-Manin \`a la mani\`ere de Greuel 
\cite{Gr} pour  $ \omega_i \in V.$ Nous notons d'ailleurs par $\psi_i$ une 
$n-$ 
forme holomorphe sur $X \setminus Cr(f)$ telle que 
$$ df_1 \wedge \cdots \wedge df_k  \wedge \psi_i = \phi_i(x) dx, 
\;\;\; 1 \leq i \leq 
\mu,$$ pour une base $\phi_j(x)dx\in 
\Blat/(f_1, \cdots, f_k) \cong \Phi.$

Alors on peut d\'eduire de la Proposition ~\ref{prop1} la relation suivante: 
$$ d{\omega}_j =  {d i_E  (\tilde{\omega}_j)} 
\leqno(1.6.1)$$ 
$$\equiv (\sum_{q=1}^{\mu} P_{jq}^{(1)} df_1 \wedge \psi_q 
-\sum_{q=1}^{\mu}  P_{jq}^{(2)} 
df_2 \wedge \psi_q + \cdots + (-1)^{k-1}\sum_{q=1}^{\mu}  
P_{jq}^{(k)} df_k \wedge \psi_q  )\; mod ((df_1,\cdots, df_k) d\Omega_X^{n-1}).
$$ 
La relation (1.6.1) implique que la d\'eriv\'ee de la forme $\omega_j$
satisfait, 
$$ d{\omega}_j \equiv
 df_1 \wedge \beta_j^{(1)} +  df_2 \wedge \beta_j^{(2)}
+ \cdots +  df_k \wedge \beta_j^{(k)} \; mod((df_1,\cdots, df_k)
d\Omega_X^{n-1}),  \leqno(1.6.2)$$ 
avec des formes m\'eromorphes $ \beta_j^{(i)}$ qui poss\`edent 
leurs p\^oles le long du lieu critique 
$Cr(f).$
 
La relation  (1.6.2) est une expression du syst\`eme de Gauss-Manin \`a la 
Greuel p.249 \cite{Gr} adopt\'ee \`a notre situation. Pour le voir, 
on remarque:
$$ \beta_j^{(i)} \equiv (-1)^{j-1}\ell_j [\sum_{q=1}^{\mu} P^{(i)}_{jq}\psi_q]
\;\;\;mod (d\Omega_X^{n-1}).$$ 
L'\'enonc\'e sur les 
p\^oles des formes $\beta_j^{(i)}$ d\'ecoule  du fait que 
$ df_1 \wedge \cdots \wedge  df_k \wedge \psi_q  \in 
\Omega_{X}^m.$ 
Voir le lemme 1.12 et la discussion \`a la p.249 de \cite{Gr}. 
 
{\bf 1.7} 
D\`es que l'expression (1.6.1) ne donne que la relation entre $d\omega_j$ 
et $\psi_j,$ elle est peu convenable pour le calcul concret du syst\`eme de 
Gauss-Manin. Il est donc 
souhaitable d'\'etablir la relation entre $d\psi_j$ et $\psi_j,$ ou bien 
$d\omega_j$ et $\omega_j.$ Dans ce but, on va chercher des relations entre 
 $\omega_j$ et $\psi_j.$ 
Si on applique $i_E$ du cot\'e gauche au (1.6.1),
$$ i_E(d\omega_j) = i_E di_E(\tilde{\omega_j})= 
(i_Ed + di_E)i_E(\tilde{\omega}_j)= \ell_j \omega_j \leqno(1.7.1)$$
$$
= \ell_j i_E(\tilde{\omega}_j) \equiv \sum_{i=1}^k (-1)^{i-1}[p_i 
\sum_{q=1}^{\mu} P_{jq}^{(i)} f_i \psi_q 
-\sum_{q=1}^{\mu} P_{jq}^{(i)} df_i \wedge i_E(\psi_q )]$$ 
$$mod((df_1,\cdots, df_k)\Omega_X^{n-1}, ((f_1,\cdots, f_k)
d\Omega_X^{n-1}).$$
Ici on a utilis\'e la formule $i_E(i_E(\omega))=0$.
  
{\bf 1.8}
La situation ci-dessus se simplifie si 
l'on regarde la  
relation 
entre des int\'egrales $\int_{\gamma(s)}\psi_q, $ au lieu de celle entre 
des formes. 
On d\'efinit l'int\'egrale- fibre $I_{\phi_q, \gamma}$ 
prise le long d'un 
cycle \'evanescent $\gamma$ dont l'ambigu\"it\'e dans l'homologie $H_n(X_s)$ 
ne sera precis\'ee qu'ult\'erieurement (voir \S 4 Th\'eor\`eme 
~\ref{thm43}), 
$$ I_{\phi_q, \gamma}(s): = \int_{\gamma(s)}\psi_q 
= (\frac{1}{2\pi i})^k  \int_{\partial \gamma(s)}\frac 
{df_1 \wedge \cdots \wedge df_k  \wedge \psi_q}{(f_1 -s_1) \cdots (f_k - s_k)}$$ 
$$ = (\frac{1}{2\pi i})^k  \int_{\partial \gamma(s)}\frac 
{\phi_q dx }{(f_1 -s_1)\cdots (f_k - s_k)}, \leqno(1.8.1)$$ 
o\`u ${\partial \gamma(s)}$ $\in H_n(X \setminus X_s)$ est un cycle 
obtenu \`a l'aide de $\partial,$ 
l'op\'erateur de cobord de Leray. Quant \`a l'op\'eration de Leray, 
on renvoie au livre de F.Pham \cite{Ph}, ou bien \`a celui de  V.A.Vasiliev 
\cite{Vas}.
 
{\bf 1.9}
De (1.7.1) on d\'eduit:
$$\ell_j \int_{\gamma(s)}{\omega}_j = 
\sum_{q=1}^{\mu} [\sum_{i=1}^k (-1)^{i-1} p_i  s_i 
P_{jq}^{(i)}(s) ] I_{\phi_q,\gamma}(s).    \leqno(1.9.1)$$ 
Cette relation est une cons\'equence imm\'ediate 
d'application 
de la d\'efinition de l'int\'egrale- fibre (1.8.1) \`a (1.7.1): 
$$ 
\int_{\gamma(s)}  df_1 \wedge i_E(\psi_j) 
=(\frac{1}{2\pi i})^k \int_{\partial \gamma(s)}\frac 
{df_1 \wedge\cdots \wedge df_k}{(f_1 -s_1)\cdots (f_k - s_k)}\wedge df_1 
\wedge i_E(\psi_j) =0.$$ 
D'une fa\c{c}on analogue, $$\int_{\partial\gamma(s)} 
\frac{df_1 \wedge\cdots \wedge df_k}
{(f_1 -s_1)\cdots (f_k - s_k)} df_i \wedge i_E(\psi_j)  =0, 
\;2 \leq i \leq k,$$
$$\int_{\partial\gamma(s)} 
\frac{df_1 \wedge\cdots \wedge df_k}
{(f_1 -s_1)\cdots (f_k - s_k)}  \wedge (f_j d\omega)  
=0, 
\;1 \leq j \leq k, \; \omega\in \Omega^{n-1}_X.$$ 
On fait comparaison entre la relation 
$$ d\int_{\gamma(s)}{\omega}_j = 
\sum_{q=1}^{\mu} [\sum_{i=1}^k (-1)^{i-1} P_{jq}^{(i)}(s)ds_i ]
I_{\phi_q, \gamma}(s), \leqno(1.9.2)  $$ 
obtenue de (1.6.1) avec la relation (1.9.1). Pour d\'eriver
(1.9.2) de (1.6.1), on utilise l'\'egalit\'e:		
$$\frac{\partial}{\partial s_\ell}\int_{\gamma(s)}{\omega}_j
= \int_{\gamma(s)}\frac{d{\omega}_j}{df_\ell},\;\; \ell= 
1,\cdots, k.$$
Voir \cite{Gr}.

En r\'esultat nous obtenons les \'equations suivantes entre 
$I_{\phi_q}(s)$ et $\frac{\partial}{\partial s_\ell} I_{\phi_q},$ $1 \leq \ell 
\leq k$ (on se passe de pr\'eciser le cycle $\gamma(s)$ sinon des cas 
exig\'es):  
$$ \sum_{q=1}^{\mu} [\sum_{i=1}^k (-1)^{i-1}p_i s_i P_{jq}^{(i)} ]
\frac{\partial}{\partial s_\ell } I_{\phi_q } \leqno(1.9.3)$$
$$= \sum_{q=1}^{\mu}  
((\ell_j -p_\ell)P_{jq}^{(\ell)}  -  p_\ell  s_\ell 
\frac{\partial}{\partial s_\ell}P_{jq}^{(\ell)} + \sum_{i\not = \ell}
(-1)^{i-1} p_i s_i 
\frac{\partial}{\partial 
s_\ell} P_{jq}^{(i)} ) I_{\phi_q },\;1\leq  j  \leq 
\mu.$$ 
C'est un syst\`eme d'\'equations qui donnent la connexion (syst\`eme) 
de Gauss-Manin. 

{\bf 1.10} 
Pour \'enoncer la proposition dans une forme plus simple, nous introduisons les 
notations suivantes: $ {\bf I}_V = $ 
$(\int{\omega}_1,$ 
$ \cdots, \int{\omega}_\mu),$ $ {\bf I}_{\Phi} =$ 
$(I_{\phi_1}(s), \cdots,$ $I_{\phi_{\mu}}(s)).$ 
On introduit les $\mu \times \mu$ matrices d\'efinies comme suit: 
$$ L_V = diag (\ell_1, \cdots, \ell_{\mu})$$ avec  $\ell_i  =  w 
({\omega}_i),$ $P^{(1)}(s) = (P^{(1)}_{jq}(s)), \cdots,$ $P^{(k)}(s)
= (P^{(k)}_{jq}(s)),
1\leq j,q  \leq \mu.  $ 
 
\begin{thm} 
 
1. 
Pour une application quasihomog\`ene 
$$ f: X \rightarrow S$$ 
aux singularit\'es isol\'ees 
d'intersection compl\`ete de dimension $n,$ le syst\`eme de 
Gauss-Manin pour 
${\bf I}_{\Phi}$ est d\'ecrit 
par les syst\`emes suivants: 
$$d[\sum_{i=1}^k (-1)^{i-1} p_i s_i P^{(i)}(s){\bf I}_{\Phi}] = 
L_V[\sum_{i=1}^k (-1)^{i-1} P^{(i)}(s)ds_i] 
{\bf I}_{\Phi}, \leqno(1.10.1)$$
ou bien, 
$$ (\sum_{i=1}^k (-1)^{i-1} p_i s_i P^{(i)}(s))\frac{\partial}{\partial 
s_\ell} {\bf I}_{\Phi} =  [L_V P^{(\ell)}(s)   -  
\frac{\partial}{\partial s_\ell}(\sum_{i=1}^k (-1)^{i-1} p_i s_i P^{(i)}(s))
] {\bf I}_{\Phi}, \leqno(1.10.2)$$ 
$1 \leq \ell \leq k.$
 
2.  La valeur critique $D$ de  d\'eformation  $X_s$  est 
donn\'e par $D=\{s \in S: \Delta(s) =0\}$ avec 
$$ \Delta(s) = det (\sum_{i=1}^k (-1)^{i-1} p_i s_i P^{(i)}(s)). \leqno(1.10.3)$$ 
 
3. Le syst\`eme (1.10.1) est un syst\`eme holon\^ome 
d'\'equations diff\'erentielles.
\label{thm03} 
\end{thm}

{\bf D\'emonstration} 

1. Dans l'expression introduite, la d\'emarche not\'ee ci-dessus peut \^etre 
interpr\'et\'ee comme suit. La relation (1.9.2) signifie 
$$ d{\bf I}_V = (\sum_{i=1}^k (-1)^{i-1} P^{(i)}(s)ds_i){\bf I}_{\Phi}. 
\leqno(1.10.4)$$ 
En revanche la relation (1.9.1) entra\^{\i}ne  
$$ L_V {\bf I}_V = (\sum_{i=1}^k (-1)^{i-1} p_i s_i P^{(i)}(s))
{\bf I}_{\Phi}. \leqno(1.10.5)$$ 
En prenant la d\'eriv\'ee de (1.10.5) et comparant celle-ci avec (1.10.3), 
on obtient la relation entre ${\bf I}_{\Phi}$ et $d{\bf I}_{\Phi}$ 
dont on peut d\'eduire (1.10.1) et (1.10.2). 
 
2. 
Il est \'etabli par Greuel que le syst\`eme de Gauss-Manin associ\'e \`a
$X_s$ poss\`ede 
son p\^ole le long de la valeur critique  de l'application $f.$ 
D'autre part, il est 
clair que (1.10.1) et (1.10.2) s'\'ecrivent comme des syst\`emes de Pfaff 
avec le p\^ole $D= \{ s \in {\bf C}^k ; det (\sum_{i=1}^k (-1)^{i-1} p_i s_i 
P^{(i)}(s)) =0 \}.$ 

3. Des \'enonnc\'es ci-dessus, il est \'evident que la vari\'et\'e 
caract\'eristique de l'\'equation (1.10.1) est un fibr\'e cotangent:
$$ T^{\ast}_D S = \{(s,\sigma) \in T^{\ast}S; \Delta(s)=0, 
<\sigma, grad\;\Delta(s)>=0\}.$$
Puisque $dim T^{\ast}_D S =k,$ (1.10.1) est un syst\`eme holon\^ome avec une 
vari\'et\'e caract\'eristique lagrangienne. 
{\bf C.Q.F.D.}

Il  faut  remarquer  ici  que  le  syst\`eme   de   Gauss-Manin   est 
compl\`etement  d\'etermin\'e  par  les  matrices  $P^{(i)}(s),$ 
$1 \leq i \leq k$
introduites dans la Proposition ~\ref{prop1}.

\vspace{2pc} 
{ 
\center{\section{Liste des syst\`emes de Gauss-Manin pour les 
singularit\'es isol\'ees simples 
 d'intersection compl\`ete de courbe espace
}} 
} 
 
Dans cette section,  on calcule le syst\`eme de Gauss-Manin associ\'e aux 
singularit\'es isol\'ees simples 
 d'intersection compl\`ete, dans le cas important celui de la
  courbe espace i.e. $n=1, k=2, m=3.$  La forme normale des 
singularit\'es isol\'ees 
simples 
 d'intersection compl\`ete (SISIC) a \'et\'e obtenue 
par M.Giusti \cite{Gui}.
Par la suite, on \'etablit une liste des notions n\'ec\'essaires pour 
d\'ecrire   le  syst\`eme   de   Gauss-Manin  associ\'e aux SISIC comme   
(1.10.1) et (1.10.2). 
 
0. Polyn\^omes $f_1$ et $f_2,$ 
 
1.  Poids  des  variables, $  w_1:=w(x_1),w_2 := w(x_2),w_3:= w(x_3), 
p_1:= w(f_1), 
p_2:=w(f_2), $ 
 
2.  L'espace vectoriel $F,$ 
 
3. L'espace vectoriel  $\tilde \Phi$ d\'efini dans $\bf 1.3$, 
 
4. Les  matrices  $P^{(1)}$  et  $P^{(2)},$ 
 
5. La fonction d\'efinissant la valeur critique 
$\Delta(s)$ pour la  d\'eformation $X_s$.

Pour la 
description de la matrice $P^{(2)}(s)$ du 4. ci-dessus, on se sert 
d'une expression comme suit: 
$$   P^{(2)}(s)= Q(s) \times V \times  \Sigma, \leqno(2.1)$$ 
o\`u les matrices composantes sont dans $GL(\mu, {\bf C}[s]).$ 
Notamment, 
$Q(s)$ indique une matrice diagonale avec 
les \'el\'ements monomiaux en les variables $s_1,s_2;$ 
 $V$ est une matrice diagonale d'\'el\'ements rationnels, 
$ \Sigma = 
\left [ 
\begin{array}{cccc} 
\sigma_{1}&0&0&0\\ 
0&\sigma_{2}&0&0\\ 
0&0&\ddots&0\\ 
0&0 &0&\sigma_{m} \\ 
\end{array} 
\right].$ 
 
Ici on a not\'e par $\sigma_{k}\in SL(\nu_k, {\bf Z}), 1 \leq k \leq m$ 
la matrice de 
permutation 
d'ordre ${\nu_k},$ 

$$ \sigma_{k} = \left [ 
\begin{array}{ccccc} 
0&0&\cdots&0&1\\ 
1&0&\cdots&0& 0\\ 
0&1&\cdots&\vdots&\vdots\\
\vdots&\vdots&\ddots&0& 0  \\ 
0 &0&\cdots &1&0 \\ 
\end{array} 
\right],$$ 
telle que $   \sigma_{k}^{\nu_k}= id_{\nu_k}. $ 
Dans la suite, on d\'ecrit  $P^{(2)}(s)$ par les donn\'ees $Q(s), V$ et $\nu_1,\nu_2 , \cdots, \nu_m$ telles que
$\sum_{i=1}^{m} \nu_k = \mu.$
 Dans les cas ci-dessous, les matrices $P^{(1)}$ et $P^{(2)}$ sont 
toutes les deux matrices semblables \`a des matrices 
diagonales.
On a choisi la num\'erotation de la base de $\tilde \Phi$ de 
sorte que $P^{(1)}$ soit une matrice diagonale. Par la suite, nous notons 
tout simplement $dx_1dx_2, dx_2dx_3$ etc. au lieu de $dx_1 \wedge dx_2, dx_2
\wedge dx_3$ afin d'\'economiser les colonnes.

 {\bf Le cas $S_{2m+3}, m\geq 1.$}

\begin{flushleft} 
 
0. \[ 
\left\{ 
\begin{array}{ccccc} 
f_1(x_1,x_2,x_3) & = & x_1^2 + x_2^2 +x_3^{2m} &=  & 0\\ 
f_2(x_1,x_2,x_3) & = & x_2x_3 & = &0 \\ 
\end{array} 
\right. 
\]

1.$$  w(x_1)=m, w(x_2)=m, w(x_3)=1, p_1:= w(f_1)= 2m, p_2:=w(f_2)=m+1. $$ 
2. $$F  = \{ x_3 
dx_1  dx_2 ,\overbrace{ x_3   dx_3   dx_1, x_3^3 dx_3   dx_1,\cdots, x_3^{2m-1}dx_3   dx_1}^{m},  dx_1  dx_2,$$
$$  
\overbrace{dx_3   dx_1, x_3^2dx_3   dx_1, \cdots,
x_3^{2m-2}dx_3   dx_1}^{m}, dx_2  dx_3 \}.$$ 
3. $$\tilde \Phi = \{\overbrace{ 1,x_3^2, \cdots, x_3^{2m}}^{m+1}, x_2, \overbrace{x_3,  x_3^3, \cdots, x_3^{2m-1}}^m, x_1 \}.$$ 
$$w(\phi_j)= \{ \overbrace{0,2, \cdots, 2m}^{m+1},m,\overbrace{1,3, \cdots, 2m-1}^{m},m  \}.$$ 
4. $ P^{(1)} = diag(s_2,\overbrace{1, 1,\cdots, 1}^{2m+1},0)$ 
 
$Q= diag(1,\overbrace{s_2,\cdots, s_2}^{m},1,1
\overbrace{s_2,\cdots, s_2}^{m-1},1), \nu_1=\nu_2=m+1, \nu_3=1.$ 
 $V = diag(2m,\overbrace{2,\cdots, 2}^{m},2m,
\overbrace{2,\cdots, 2}^{m+1})$
 
5. $$ \Delta(s)=s_2^2 ((\frac{-s_1}{m+1})^{m+1} - (\frac{s_2^2}{m})^{m})^2. $$ 
\end{flushleft}

 {\bf Le cas $S_{2m+4}, m \geq 1$}
 
\begin{flushleft} 
 
0. \[ 
\left\{ 
\begin{array}{ccccc} 
f_1(x_1,x_2,x_3) & = & x_1^2 + x_2^2 +x_3^{2m+1} &=  & 0\\ 
f_2(x_1,x_2,x_3) & = & x_2x_3 & = &0 \\ 
\end{array} 
\right. 
\]

1.$$  w(x_1)=2m+1, w(x_2)=2m+1, w(x_3)=2, p_1:= w(f_1)= 2(2m+1), p_2:=w(f_2)= 2m+3. $$ 
2. $$F = \{ x_3 
dx_1  dx_2 ,\overbrace{ x_3   dx_3   dx_1, x_3^3 dx_3   dx_1,\cdots, x_3^{2m-1}dx_3   dx_1}^{m},$$ $$  dx_1  dx_2,  
\overbrace{dx_3   dx_1, x_3^2dx_3   dx_1, \cdots,
x_3^{2m}dx_3   dx_1}^{m+1}, dx_2  dx_3 \}.$$ 
3. $$\tilde \Phi = \{\overbrace{ 1, x_3^2, \cdots, x_3^{2m}}^{m+1}, x_2, 
\overbrace{x_3,  x_3^3, \cdots, x_3^{2m+1}}^{m+1}, x_1 \}.$$ 
$$w(\phi_j)= \{ \overbrace{0,4, \cdots, 4m}^{m+1},
2m+1 ,\overbrace{2,6 , \cdots, 2(2m+1)}^{m+1},
2m+1  \}.$$
4. $ P^{(1)} = diag(s_2,\overbrace{1, 1,\cdots, 1}^{2m+2},0)$ 
 
$Q= diag(1,\overbrace{s_2,\cdots, s_2}^{m},1,1,
\overbrace{s_2,\cdots, s_2}^{m},1), \nu_1=2m+3, \nu_2=1.$ 
 $V= diag(2m+1,\overbrace{2,\cdots, 2}^{m}, 2m+1,\overbrace{2,\cdots, 2}^{m +2} ).$

5. $$ \Delta(s)
=s_2^2 ((\frac{s_1}{2m+3})^{2m+3}+ (\frac{s_2^2}{2m+1})^{2m+1}). 
$$ 
\end{flushleft}

{\bf Le cas $T_7.$}

\begin{flushleft} 
 
0. \[ 
\left\{ 
\begin{array}{ccccc} 
f_1(x_1,x_2,x_3) & = & x_1^2 + x_2^3 +x_3^3 &=  & 0\\ 
f_2(x_1,x_2,x_3) & = & x_2x_3 & = &0 \\ 
\end{array} 
\right. 
\] 
  
1.$$  w(x_1)=3, w(x_2)=2, w(x_3)=2, p_1:= w(f_1)= 6, p_2:=w(f_2)=4. $$ 
2. $$F  = \{ x_3 
dx_1  dx_2 , x_3^2 dx_3   dx_1, dx_1  dx_2, x_3 dx_3   dx_1,  x_2 
dx_1  dx_2,  dx_3   dx_1,  dx_2  dx_3 \}.$$ 
3. $$\tilde \Phi = \{ 1, x_3^3, x_2, x_3^2, x_2^2, x_3, x_1 \}.$$ 
$$w(\phi_j)=  \{ 0,6,2,4,4,2,3 \}.$$ 
4. $ P^{(1)} = diag( s_2, 1,1,1,1,1,0)$ 
 
$Q= diag(s_2^2,1,1,s_2,1,s_2,1,s_2,1), \nu_1=\nu_2= \nu_3=2, \nu_4=1.$ 
$V= diag(3,3,3,3,3,3,2)$
 
5. 
$$ \Delta(s)=s_2^2 (s_1^2 - 4s_2^3)^3. $$ 
\end{flushleft}

{\bf Le cas $T_8.$} 
\begin{flushleft} 
 
0.\[ 
\left\{ 
\begin{array}{ccccc} 
f_1(x_1,x_2,x_3) & = & x_1^2 + x_2^3 +x_3^4 &=  & 0\\ 
f_2(x_1,x_2,x_3) & = & x_2x_3 & = &0 \\ 
\end{array} 
\right. 
\] 
 
1. 
$$  w(x_1)=6, w(x_2)=4, w(x_3)=3, p_1:= w(f_1)= 12,p_2:=w(f_2)=7.$$ 
2. 
$$F = \{ x_3 dx_1  dx_2 , x_3^2 dx_3   dx_1, dx_1  dx_2, x_3  dx_3   dx_1,  x_2 
dx_1  dx_2, dx_3   dx_1, x_3^3 dx_3   dx_1,  dx_2  dx_3 \}.$$ 
3. 
$$\tilde \Phi = \{ 1, x_3^3, x_2, x_3^2, x_2^2, x_3, x_3^4, x_1 \}.$$ 
$$w(\phi_j)=\{ 0,9,4,6,8,3,12,6 \}.$$ 
4. 
$$ P^{(1)} = diag( s_2, 1,1,1,1,1,1,0)$$ 
  
$Q= diag(1, s_2^2,1,s_2,s_2,1,s_2^2,1), \nu_1=7, \nu_2=1.$  
$V= diag(4, 3,4,3,4,3,3,1)$

5. 
$$ \Delta(s)=s_2^2 (3^3 4^4 s_1^7 - 7^7 s_2^{12}). $$ 
\end{flushleft} 
 
{\bf Le cas $T_9.$} 
\begin{flushleft} 
0.\[ 
\left\{ 
\begin{array}{ccccc} 
f_1(x_1,x_2,x_3) & = & x_1^2 + x_2^3 +x_3^5& =  & 0\\ 
f_2(x_1,x_2,x_3) & = & x_2x_3 & = &0 \\ 
\end{array} 
\right. 
\]

1. 
$$  w(x_1)=15, w(x_2)=10, w(x_3)=6, p_1:= w(f_1)= 30, p_2:=w(f_2)=16 . $$ 
2. 
$$F = \{ x_3 dx_1  dx_2 ,  x_3^2  dx_3   dx_1,  x_2  dx_1  dx_2,  dx_3   dx_1, $$
$$
x_3^3dx_3   dx_1,  dx_1  dx_2, x_3   dx_3   dx_1, x_3^4 dx_3   dx_1,  dx_2  dx_3 
\}.$$ 
3. 
$$\tilde \Phi = \{ 1, x_3^3, x_2^2, x_3,  x_3^4, x_2, x_3^2, x_3^5, x_1 \}.$$ 
$$w(\phi_j)=\{ 0,18,20,6,24,10,12,30, 15 \}.$$ 
4. 
$$ P^{(1)} = diag( s_2, 1,1,1,1,1,1,1,0)$$ 
$Q= diag(1, s_2^2, s_2,1, s_2,1,s_2,s_2^2,1), \nu_1=8, \nu_2=1.$ 
$V = diag(5, 3,5,3, 3,5,3,3,2)$

5. 
$$ \Delta(s)=s_2 (5^5 3^3s_1^8 - 2^{24} s_2^{15}). $$ 
\end{flushleft}

 {\bf Le cas $U_7.$} 
\begin{flushleft} 
0. 
\[ 
\left\{ 
\begin{array}{ccccc} 
f_1(x_1,x_2,x_3) & = & x_1^2 + x_2x_3 &=  & 0\\ 
f_2(x_1,x_2,x_3) & = & x_1x_2 + x_3^3 & = &0 \\ 
\end{array} 
\right. 
\] 
 
1. 
$$  w(x_1)=4, w(x_2)=5, w(x_3)=3, p_1:= w(f_1)=8,p_2:=w(f_2)=9.$$ 
2. 
$$F = \{ {x_3 dx_2  dx_3 + x_1 dx_3 dx_1 } , x_3^2 dx_3 dx_1 + x_1 
dx_1  dx_2, $$  $$x_3^2 dx_2  dx_3 + x_1x_3 dx_3 dx_1, dx_3 dx_1,  dx_2  dx_3, dx_1  dx_2, x_3 dx_3 dx_1 
\}.$$ 
3. 
$$\tilde \Phi = \{ 1, x_1x_3^2,x_3, x_1, x_2,  x_3^2, x_1x_3 \}.$$ 
$$w(\phi_j)=\{ 0,10,3,4,5,6,7 \}.$$ 
4. 
$$ P^{(1)} = diag( s_1, 1,s_1,1,1,1,1)$$ 
$Q= diag(1, s_2, 1, 1,1 ,1,1), \nu_1=7.$ 
$V = diag(3, \frac{1}{4}, 3, 1,2,\frac{1}{3},1)$

5. 
$$ \Delta(s)=2^{22} s_1^9 - 3^{15} s_2^8. $$ 
\end{flushleft} 
 
{\bf Le cas $U_8.$} 
 
\begin{flushleft} 
 
0.\[ 
\left\{ 
\begin{array}{ccccc} 
f_1(x_1,x_2,x_3) & = & x_1^2 + x_2x_3 + x_3^3&=  & 0\\ 
f_2(x_1,x_2,x_3) & = & x_1x_2 & = &0 \\\ 
\end{array} 
\right. 
\] 
 
1. 
$$  w(x_1)=3, w(x_2)=4, w(x_3)=2, p_1:= w(f_1)=6,p_2:=w(f_2)=7. $$ 
2. 
$$F = \{ x_1 dx_2  dx_3 , x_1x_3 dx_3 dx_1 -\frac{x_1}{3}dx_1  dx_2, $$ 
  $$x_1x_3 dx_2  dx_3,  dx_3 dx_1,  dx_2  dx_3,x_3 dx_3 dx_1 -\frac{1}{3}dx_1  dx_2,  x_1  dx_3 dx_1, dx_1  dx_2 
\}.$$ 
3. 
$$\tilde \Phi = \{ 1, x_1^2x_3, x_3, x_1, x_2, x_1x_3, x_1^2, x_3^2 \}.$$ 
$$w(\phi_j)=\{ 0,8,2,3,4,5,6,4 \}.$$ 
4. 
$$ P^{(1)} = diag( s_2, 1,s_2,1,1,1,1,0)$$ 
$Q= diag(1, s_2, 1,1,1,1,1,1), \nu_1=7, \nu_2=1.$
$V = diag(2, \frac{-1}{3}, 2, 1,2,\frac{-1}{3},1,3)$

5. 
$$ \Delta(s)= s_2^3(2^4 3^9 s_1^7 -  7^7 s_2^6). $$ 
 
\end{flushleft} 

{\bf Le cas $U_9.$} 
\begin{flushleft} 
0.\[ 
\left\{ 
\begin{array}{ccccc} 
f_1(x_1,x_2,x_3) & = & x_1^2 + x_2x_3 &=  & 0\\ 
f_2(x_1,x_2,x_3) & = & x_1x_2 + x_3^4& = &0 \\\ 
\end{array} 
\right. 
\] 
 
1. 
$$  w(x_1)=5, w(x_2)=7, w(x_3)=3, p_1:= w(f_1)= 10,p_2:=w(f_2)=12. $$ 
2. 
$$F= \{  {x_3 dx_2  dx_3 } + x_1 dx_3   dx_1 , 
x_3^3 dx_3 dx_1 +x_1 dx_1  dx_2, \ 
x_3^3 dx_2  dx_3 + x_1x_3^2 dx_3   dx_1, $$ $$x_3 dx_3 dx_1, 
dx_2  dx_3, dx_1  dx_2, x_3^2 dx_3 dx_1, x_3^2 dx_2  dx_3 + x_1x_3 dx_3   dx_1,  dx_3 dx_1, 
\}.$$ 
3. 
$$\tilde \Phi = \{ 1, x_1x_3^3,x_3^2, x_1x_3, x_2, x_3^3,x_1x_3^2, x_3,  x_1 \}.$$ 
$$w(\phi_j)=\{ 0,14,6,8,7,9,11,3,5 \}.$$ 
4. 
$$ P^{(1)} = diag( s_1, 1,s_1,1,1,1,s_1,1)$$ 
$Q= diag(1, s_2, 1, 1, 1, 1,1,1,1), \nu_1=4, \nu_2=5.$
$V= diag(3, \frac{1}{5}, 3, 1, 8,\frac{1}{4},1,3, 1)$

5. 
$$ \Delta(s)= (5^5 s_1^6 - 2^4 3^6 s_2^5)^2. $$ 
\end{flushleft}

{\bf Le cas $W_8.$} 
 
0. 
\begin{flushleft} 
\[ 
\left\{ 
\begin{array}{ccccc} 
f_1(x_1,x_2,x_3) & = & x_2^2 + x_1x_3  & =  & 0\\ 
f_2(x_1,x_2,x_3) & = & x_1^2 + x_3^3 & = &0 \\ 
\end{array} 
\right. 
\] 
 
1. 
$$  w(x_1)= 6, w(x_2)= 5, w(x_3)= 4, p_1:= w(f_1)= 10, 
p_2:=w(f_2)= 12. $$ 
 
2. 
$$F = \{ 3x_1 dx_2  dx_3 + 2x_3 dx_1  dx_2,  x_3^2 dx_2  dx_3 + x_1 
dx_1  dx_2,$$ $$ 
3x_1x_3 dx_2  dx_3 + 2x_3^2 dx_1  dx_2, 
dx_2  dx_3, dx_1  dx_2, x_3   dx_2  dx_3, dx_3   dx_1, x_3   dx_3   dx_1 \}.$$ 
3. 
$$\tilde \Phi = \{ 1, x_1x_3^2,x_3, x_1, x_3^2,x_1x_3,  x_2, x_2x_3 \}.$$ 
$$w(\phi_j)=\{ 0,14,4,6,8,10,5,9 \}.$$ 
4. 
$$ P^{(1)} = diag( 6s_2, 5, 6s_2,1,1,2,0,0)$$ 
$Q= diag(1, s_2, 1,1,1,1,1,1), \nu_1=6, \nu_2= 
\nu_3 =1.$
$V= diag(5, 1, 5, \frac{1}{2},\frac{1}{3},1,1,2)$

5. 
$$ \Delta(s)=s_2^4 (5^5 s_1^6 -  2^2 3^3 s_2^{5}). $$

\end{flushleft} 
 
{\bf Le cas $W_9.$} 
 
0.\begin{flushleft} 
\[ 
\left\{ 
\begin{array}{ccccc} 
f_1(x_1,x_2,x_3) & = & x_2^2 + x_1x_3  & =  & 0\\ 
f_2(x_1,x_2,x_3) & = & x_1^2  + x_2x_3^2  & = &0 \\ 
\end{array} 
\right. 
\] 
 
1. 
$$  w(x_1)= 5, w(x_2)= 4, w(x_3)= 3, p_1:= w(f_1)= 8, p_2:=
w(f_2)=10. $$ 
2. 
$$F = \{x_1 dx_2  dx_3 + x_3 dx_1  dx_2, x_2x_3 dx_2  dx_3 + x_1 
dx_1  dx_2, x_1x_2 dx_2  dx_3 + x_2 x_3 dx_1  dx_2,$$  $$  dx_3   dx_1, x_3 dx_2  dx_3, 2x_1 dx_3   dx_1 + x_2 
dx_1  dx_2, dx_2  dx_3, dx_1  dx_2, x_2  dx_2  dx_3 
\}.$$ 
3. 
$$\tilde \Phi = \{ 1, x_1x_2x_3, x_2, x_3^2, x_1x_3, x_3, x_1,   x_2x_3, x_1x_2 \}.$$ 
$$w(\phi_j)=\{ 0,12,4,6,8,3,5,7,9 \}.$$ 
4. 
$$ P^{(1)} = diag( 2s_2, 4, 2s_2,1,2,s_1,2,2,2)$$ 
$Q=diag( 1, s_2, 1,\cdots, 1), \nu_1=5, \nu_2=4.$ 
 $V = diag(2,1,2,2,1,5,1,1,1)$

5. 
$$ \Delta(s)=s_2^3 (2^{12}s_1^5 -  5^5 s_2^{4})^2. $$ 
 \end{flushleft}
{\bf  Le cas $Z_9.$} 
 
\begin{flushleft} 
 
0.\[ 
\left\{ 
\begin{array}{ccccc} 
f_1(x_1,x_2,x_3) & = & x_1^2 + x_2^2 + x_3^3 & =  & 0\\ 
f_2(x_1,x_2,x_3) & = & x_1x_2  & = &0 \\ 
\end{array} 
\right. 
\] 
 
1. 
$$  w(x_1)= 3, w(x_2)= 3, w(x_3)= 2, p_1:= w(f_1)= 6, 
p_2:=w(f_2)= 6. $$ 
2. 
$$F  =  \{  x_1  dx_2  dx_3,   x_1  dx_3   dx_1,   dx_2  dx_3,   dx_3   dx_1, 
x_1x_3   dx_2  dx_3,$$ $$ x_1x_3   dx_3   dx_1, x_3 dx_2  dx_3, x_3   dx_3   dx_1, dx_1  dx_2 
\}.$$ 
3. 
$$\tilde \Phi = \{ 1, x_1^2, x_2, x_1, x_3, x_1^2x_3, x_2x_3, x_1x_3, x_3^2 \}.$$ 
$$w(\phi_j)=\{ 0,6,3,3,2,8,5,5,4 \}.$$ 
 
4. 
$$ P^{(1)} = diag( s_2, 1,1,1, s_2,1,1,1,0)$$ 
$Q= diag(1,s_2,1,1,1,s_2,1,1,1), \nu_1= \nu_2 = \nu_3 = \nu_4 =2, \nu_5=1.$ 
$V = diag(2,2,2,2,2,2,2,2,1)$

5. 
$$ \Delta(s)=s_2^3 (s_1^2 -  4 s_2^{2})^4. $$ 
\end{flushleft}

{\bf Le cas $Z_{10}.$} 
\begin{flushleft} 
0.\[ 
\left\{ 
\begin{array}{ccccc} 
f_1(x_1,x_2,x_3) & = & x_1^2 + x_2x_3^2  & =  & 0\\ 
f_2(x_1,x_2,x_3) & = & x_2^2  + x_3^3  & = &0 \\ 
\end{array} 
\right. 
\] 
 
1. 
$$  w(x_1)= 7, w(x_2)= 6, w(x_3)= 4, p_1:= w(f_1)= 14, 
p_2:=w(f_2)= 12
. $$ 
2. 
$$F  = \{ 
3x_2 dx_3   dx_1 + 2x_3 dx_1  dx_2, x_1  dx_2  dx_3  +  x_3   dx_1  dx_2, 
x_3 dx_3   dx_1, dx_1  dx_2, dx_3   dx_1,$$ $$ x_1x_2 dx_2  dx_3 + x_2 x_3 dx_1  dx_2, 
3x_2x_3 dx_3   dx_1 + 2x_3^2 dx_1  dx_2, 2x_3^2 dx_3   dx_1 + x_2 dx_1  dx_2, 
dx_2  dx_3, x_3   dx_2 dx_3 \}.$$ 
3. 
$$\tilde \Phi = \{ 1, x_3^3, x_2x_3,x_3^2, x_2, x_2x_3^3, x_3, x_2x_3^2, x_1, x_1x_3 \}.$$ 
$$w(\phi_j)=\{ 0,12,10,8,6,18,4,14,7,11 \}.$$ 
 
4. 
$$ P^{(1)} = diag( 6s_2, 3,2,3,2,3,6s_2,7,0,0)$$ 
$Q= diag(1,s_1,1,1,1,s_1,1,s_2,1,1), \nu_1= 8, \nu_2 = \nu_3 = 1.$ 
 $V= diag(7,2,1,2,1,2,7,2,1,1)$

5. 
$$ \Delta(s)=s_1^2 s_2^4 (7^7 s_1^6 -  2^8 3^3 s_2 ^{7}). $$ 
\end{flushleft} 

{ 
\center{\section{ 
Les spectres du syst\`eme de Gauss-Manin 
}} 
} 
\vspace{2pc} 
 
Pour les singularit\'es simples SISIC, on peut mettre en \'evidence
les informations 
topologiques sur la singularit\'e \`a partir des syst\`emes 
 (1.10.1) et (1.10.2). Nous formulons  ce fait comme suivant. 
 
{\bf 3.1}
\begin{thm} 
 
1. Le syst\`eme de Gauss-Manin pour ${\bf I}_{\Phi}$ associ\'e 
aux singularit\'es isol\'ees simples d'intersection 
compl\`ete de courbe espace s'\'ecrit 
sous la forme suivante: 
$$ P'^{(1)}(s_1)( s_1 id_{\mu}\frac{\partial}{\partial 
s_1}   - \frac{1}{p_1}L_{\Phi}) {\bf I}_{\Phi} = \frac{ p_2}{p_1} s_2 
P'^{(2)}(s_2)\frac{\partial}{\partial s_1}{\bf I}_{\Phi}, \leqno(3.1.1)$$ 
o\`u $P'^{(1)}(s_1),$ une matrice diagonale, et $P'^{(2)}(s_2),$ une matrice 
semblable \`a une matrice diagonale. 
$$ \frac{1}{p_1}L_{\Phi} = diag\{\lambda_1, \cdots, \lambda_{\mu}   \}, 
\lambda_j  = \frac {w(\psi_j)}{p_1}.$$ 
D'une fa\c{c}on analogue: 
$$P'^{(2)}(s_2) ( s_2 id_{\mu}\frac{\partial}{\partial 
s_2}   - \frac{1}{p_2}L_{\Phi}) {\bf I}_{\Phi} = \frac{ p_1}{p_2} s_1 
P'^{(1)}(s_1)\frac{\partial}{\partial s_2}{\bf I}_{\Phi}. \leqno(3.1.2)$$ 
 
2.(la sym\'etrie des spectres).  Notons $\sigma \in {\cal S}_{\mu}$  
la permutation telle que: 
$$\tilde{\lambda}_1 
 \leq \tilde{\lambda}_2 \leq 
\cdots \leq \tilde{\lambda}_{\mu},$$
alors il existe un nombre rationnel $\tilde 
{\lambda}_0$ 
tel que 
$$ \tilde{\lambda}_0 - \tilde{\lambda}_i =  \tilde{\lambda}_{\mu -i}  - 
\tilde{ \lambda}_{0},\;\;\;\; 
1 \leq i \leq \mu.   $$ 
\label{thm1} 
\end{thm}

\begin{dfn}
 
Nous appelons les rationnels $ \{\tilde{\lambda_1}, 
\cdots,  \tilde{\lambda_{\mu}}  \} \in \frac{1}{p_1}{\bf Z}_{\geq 0}$ 
les spectres du syst\`eme de Gauss-Manin (3.1.1). D'une fa\c{c}on analogue, 
$ \{\frac{p_1}{p_2}\tilde{\lambda_1}, 
\cdots,  \frac{p_1}{p_2}\tilde{\lambda_{\mu}}  \} \in 
\frac{1}{p_2}{\bf Z}_{\geq 0}$ les spectres du syst\`eme  (3.1.2). 
\label{dfn31}
\end{dfn}
 
\begin{remark} 
{\em La propri\'et\'e de sym\'etrie des spectres de la  structure  de  Hodge 
mixte de la cohomologie relative associ\'ee \`a SIIC 
a \'et\'e d\'emontr\'ee par W.Ebeling et J.Steenbrink 
\cite{ESt}. Les calculs  concrets  ont  \'et\'e  achev\'es  par  ce 
dernier pour   les   singularit\'es    isol\'ees unimodales 
d'intersection compl\`ete \cite{St}.

Notre approche est diff\'erent de celui de Ebeling-Steenbrink,
puisque nos objets principaux sur lesquel la transformation de monodromie
agit sont les espaces $V$ et $\Phi$ de Greuel-Hamm. En g\'en\'eral la dimension
de la cohomologie relative est plus grande que le nombre de Milnor
de la singularit\'e $X_0$ car une structure suppl\'ementaire 
intervient dans la cohomologie relative. Notamment ils regardent 
une d\'eformation
d'une SIIC d\'ependant de deux param\`etres:
$$ (f,g):(X_s,x) \rightarrow ({\bf C}^2, 0),$$
en sorte  que  la fonction $g$ d\'efinisse une singularit\'e isol\'ee 
d'hypersurface non-d\'eg\'en\'er\'ee. C'est la cohomologie de la fibre de Milnor
de $g$ qui s'entremet dans la cohomologie relative.
 } 
\label{remark31}
\end{remark} 
 
{\bf 3.2} 
{\bf D\'emonstration} 
 
1. D'abord on observe la possibilit\'e d'\'etendre la connexion 
de Gauss-Manin sur un module plus grand que $F.$ 
On regarde un module $ F'=F[\frac{1}{f_1}, \frac{1}{f_2}]$ 
au lieu de $F,$ et choisit sa base $ \tilde{\omega_i}', 
\;(1 \leq i \leq \mu)$ de telle sorte que la relation suivante analogue 
\`a $(1.10.5)$ ait lieu pour 
$ {\bf I}_{F'} = $ 
$(\int_{\gamma(s)}i_E(\tilde{\omega}_1'),$ 
$ \cdots, \int_{\gamma(s)}i_E(\tilde{\omega}_\mu') ):$ 
$${L}_{F'}{\bf I}_{F'} = 
( p_1 s_1 P'^{(1)}(s_1) - p_2 s_2 P'^{(2)}(s_2)) {\bf I}_{\Phi}, 
\leqno(3.2.1)$$ 
o\`u $ L_{F'}= diag( \ell_1', \cdots, \ell_\mu'),$ et $ \ell_i' = 
w(\tilde{\omega}_{i}'). $
Pour le voir, on d\'efinit les 
formes 
 $ \tilde{\omega_i}', $ comme suit: 
 
\begin{center} 
$ 
\begin{array}{ccc} 
  \tilde{\omega_i}'  & = &\frac{ \tilde{\omega_i}} 
{f_1^{\eta_i} f_2^{\delta_i}} 
\;\; {\rm si} \;\frac{P^{(1)}_{i,j}(s)}{s_1^{\eta_i} s_2^{\delta_i}} \in 
{\bf C}[s_1, s_1^{-1}], {\rm et }\; 
\frac{P'^{(2)}_{i,j}(s)}{s_1^{\eta_i} s_2^{\delta_i}} \in 
{\bf C}[s_2, s_2^{-1}], \eta_i,\delta_i= 0,1,2,\cdots. \\ 
\end{array} 
$ 
\end{center} 
Nous nous servirons des notations 
$$ P'^{(1)}(s_1) = diag( s_1^{-\eta_1} s_2^{-\delta_1}, \cdots , 
s_1^{-\eta_\mu} s_2^{-\delta_\mu}) \times P^{(1)}(s_1,s_2),$$ 
$$ P'^{(2)}(s_2) = diag( s_1^{-\eta_1} s_2^{-\delta_1}, \cdots , 
s_1^{-\eta_\mu} s_2^{-\delta_\mu}) \times P^{(2)}(s_1,s_2).$$ 
C'est \`a dire:
$$ P'^{(1)}(s_1) = diag( s_1^{\tilde \eta_1}, \cdots , 
s_1^{\tilde \eta_\mu} ) \times diag( p_1^{(1)}, \cdots , 
p_\mu^{(1)} ) , \leqno(3.2.2)$$ 
$$ P'^{(2)}(s_2) = diag( s_2^{\tilde \delta_1}, \cdots , 
s_1^{\tilde \delta_\mu} ) \times diag( p_1^{(2)}, \cdots , 
p_\mu^{(2)} ) \cdot \Sigma,$$
o\`u  $p_i^{(\ell)},$ $(1\leq i \leq \mu, \ell =1,2)$ sont des 
rationnels et $\Sigma $ une matrice comme dans $(2.1).$
Cette op\'eration est faisable pour toutes les $P^{(1)}, P^{(2)}$
 calcul\'ees 
dans \S 2 car il existe au plus des entiers uniques
$j_1, j_2 \in [1,\mu]$
tels que
$$\begin{array}{c}
\tilde \omega_{j} \wedge df_1  =  P^{(2)}_{jj_2} (f) \phi_{j_2}(x) dx \\
\tilde \omega_{j} \wedge df_2  =  P^{(1)}_{jj_1} (f) \phi_{j_1}(x)dx 
\end{array}$$
pour chaque
$j \in [1,\mu].$

En bref, on arrive \`a l'expression (3.2.1), 
si on multiplie la matrice 
$ diag( s_1^{-\eta_1} s_2^{-\delta_1},$ $\cdots ,$ 
$s_1^{-\eta_\mu} s_2^{-\delta_\mu})$ du c\^ot\'e gauche \`a $(1.10.5).$ 
On obtient par une mani\`ere analogue \`a (1.6.1): 
$$\frac{\partial}{\partial s_1}  {\bf I}_{F'} = P'^{(1)}(s_1){\bf 
I}_{\Phi}. 
\leqno(3.2.3)$$

D'autre part, ayant diff\'erenti\'e l'expression (3.2.1), on obtient: 
$$\ell_i' \frac{\partial}{\partial s_1} \int_{\gamma(s)} 
i_E(\tilde{\omega_i}') 
  = \sum_{j=1}^{\mu} (p_1 P'^{(1)}_{ij} + p_1 s_1 \frac{\partial}{\partial 
s_1} P'^{(1)}_{ij}) I_{\phi_j } +  \sum_{j=1}^{\mu} (p_1 s_1  P'^{(1)}_{ij}(s_1)- 
p_2 s_2 P'^{(2)}_{ij}(s_2)) \frac{\partial}{\partial s_1}I_{\phi_j },$$ 
c'est- \`a- dire, 
$$ L_{F'}\frac{\partial}{\partial s_1}  {\bf I}_{F'} = (p_1 P'^{(1)}(s_1) 
+ p_1 s_1\frac{\partial}{\partial s_1} P'^{(1)}){\bf I}_{\Phi} + 
(p_1 s_1  P'^{(1)}(s_1)-p_2 s_2 P'^{(2)}(s_2)) \frac{\partial}{\partial s_1}{\bf 
I}_{\Phi}, 
\leqno(3.2.4)$$ 
En suite, on remarque que la relation suivante:
$$ \tilde{\omega_i'} = \sum_{1 \leq i,j \leq \mu} P'^{(1)}_{ij} \wedge df_1 
(\frac{\phi_j dx_1\wedge dx_2\wedge dx_3 }{df_1  \wedge df_2}) - 
\sum_{1 \leq i,j \leq \mu}\ P'^{(2)}_{ij}\wedge df_2 (\frac{\phi_j dx_1\wedge dx_2\wedge dx_3 
}{df_1  \wedge df_2})$$ 
donne une relation entre les \'el\'ements  de  matrices  $L_{\Phi}, 
L_{F'}$  et $P'^{(1)}:$ 
$$ \ell_i' = \frac{E_{\ast}(P'^{(1)}_{ij})}{P'^{(1)}_{ij}} + p_1 + L_{\Phi,j}
= \frac{E_{\ast}(P'^{(2)}_{ij})}{P'^{(2)}_{ij}} + p_2 + L_{\Phi,j},$$ 
ici $E_{\ast} = p_1 s_1 \frac{\partial}{\partial s_1} +p_2 s_2 
\frac{\partial}{\partial s_2}, $ le champ d'Euler sur $S$. 
Ce dernier  entra\^{\i}ne 
$$ \sum_{j=1}^{\mu}( -\ell_i' + p_1 + L_{\Phi,j}  )P'^{(\ell)}_{ij} + 
E_{\ast}(P'^{(\ell)}_{ij}) 
=0, \ell =1,2.$$ 
Autrement dit, 
$$ P'^{(1)}( L_{\Phi}- L_{F'} +  p_1\cdot id_{\mu}) + E_{\ast} (P'^{(1)}) =0 
\leqno(3.2.5)$$ 
$$ P'^{(2)}( L_{\Phi}- L_{F'} +  p_2 \cdot id_{\mu}) + E_{\ast} (P'^{(2)}) =0 
\leqno(3.2.6)$$ 
En somme, $(3.2.3), (3.2.4),$ et $ (3.2.5)$ nous m\`enent \`a conclure 
$(3.1.1).$ 
 
La d\'emonstration de (3.1.2) est parall\`elle \`a 
celle de (3.1.1), en tenant compte de $(3.2.6).$ 
 
2. 
Nous introduisons ici la notation $|w|:= w_1 + w_2 + w_3$
et $|p|: = p_1 + p_2.$ 
Puisque $\lambda_j = \frac{1}{p_1} (w(\phi_j) + |p| - |w|),$ 
il suffit de d\'emontrer la sym\'etrie entre les poids 
des \'el\'ements de l'espace $\tilde \Phi.$

On peut d\'eduire de  \cite{Al}, \cite{Al1}, 3.4 que 
la s\'erie de Poincar\'e $P_{\tilde \Phi}(t)$ 
de l'espace vectoriel $\tilde \Phi$  s'\'ecrit comme suit: 
$$    P_{\tilde \Phi}(t)    =     t^{|p|-|w|}+     
(1-t^{|p|-|w|})\frac{(1-t^{p_1})(1- 
t^{p_2})} {(1-t^{w_1})(1-t^{w_2})(1-t^{w_3})}. \leqno(3.2.6) $$ 
Il est facile de voir que 
le polyn\^ome $P_{\tilde \Phi}(t) $ 
a coefficients sym\'etriques par rapport au 
terme central  
$t^{|p|-|w|}.$ {\bf C.Q.F.D.} 
\vspace{2pc} 
{ 
\center{\section{ 
L'expression explicite de la  transform\'ee de Mellin de l'int\'egrale 
-fibre 
}} 
} 
Dans cette section, nous essayons d'\'etablir une expression explicite 
 de la  transform\'ee de Mellin de l'int\'egrale 
- fibre  au moyen des invariants topologiques de singularit\'es. 
 
{\bf 4.1. L'EDF et ses solutions explicites } 
 
D'abord on \'etablit les \'equations aux diff\'erences finies (EDF) 
pour la  transform\'ee de Mellin $M_{\phi_j}(z_1, z_2)$ de $I_{\phi_j}(s_1,s_2):$ 
$$  M_j(z_1,z_2) = \int_{\gamma} I_j(s_1,s_2) s_1^{z_1} s_2^{z_2} 
ds_1 ds_2$$ pour un certain $\gamma$ qui \'evite les p\^oles de $I_j(s_1,s_2)$ 
(on note $M_j(z)$ et $I_j(s)$   au lieu de 
$M_{\phi_j}(z_1, z_2),$ $I_{\phi_j}(s_1,s_2)$ 
pour all\'eger l'\'ecriture). Il est \'evident que l'int\'egrale $M_j(z_1,z_2)$ est bien d\'efinie 
pour $ \Re z_1, \Re z_2 >> 0, |\Im z_1|, |\Im z_2| >>0 $ 
gr\^ace \`a la  r\'egularit\'e de
tous ses points singuliers y compris l'infini
de l'\'equation diff\'erentielle satisfaite par 
les int\'egrales $I_j(s_1,s_2)$ 
(voir
Th\'eor\`eme ~\ref{thm1} ).  Cette d\'emarche est bien formul\'ee dans \cite{LS} \S1, sous le terme
de ``transformation de Mellin alg\'ebrique''.

On applique la transformation de Mellin \`a la relation (3.1.1). Alors 
on en tire l'EDF entre $M_j(z_1, z_2)$. 
A l'aide des  $\tilde{\eta}_k , \tilde{\delta}_k $ 
d\'efinis dans (3.2.2)  et $\nu_k$ introduit 
dans le \S 2, elles 
s'expriment comme suit: 
$$(z_1 +\tilde{\eta}_k + \lambda_k +1)M_k(z_1+ \tilde{\eta}_k, z_2)  
= \tilde{v}_k z_1 M_{[k-1]}(z_1-1, z_2+1 + \tilde{\delta}_k), \leqno(4.1.1)$$ 
o\`u $\tilde{v}_k \in {\bf Q}$
$ \nu_1+ \nu_2 + \cdots +\nu_m +1 \leq k \leq\nu_1+ \cdots +\nu_{m+1},  $ 
et $[k-1]= k-1$ si  $\nu_1+ \nu_2 + \cdots +\nu_m +1  < k \leq\nu_1+ 
\cdots +\nu_{m+1},$ et   
$[k-1]=\nu_1+ \nu_2 + \cdots +\nu_{m+1}$ si 
$  k =\nu_1+ \nu_2 + \cdots +\nu_m +1. $ L'\'equation (4.1.1) se d\'eduit 
du 
fait que l'int\'egrale $I_k(s_1,s_2)$ 
est li\'ee  \`a l'autre int\'egrale $I_{k'}(s_1,s_2), $ par une relation 
non-triviale si et seulement si $\nu_1+ \nu_2 + \cdots +\nu_m +1 
\leq k, k' \leq 
\nu_1+ \cdots +\nu_{m+1},$ i.e. si et seulement si elles 
sont toutes les deux d'un 
bloc de taille $\nu_{m+1}.$ Nous disons que la relation de r\'ecurrence
(4.1.1) se ferme pour $ \phi_{\nu_1+ \nu_2 + \cdots +\nu_m +1}(x),$
$\cdots,$ $\phi_{\nu_1+ \cdots +\nu_{m+1}}(x).$
  On va repr\'esenter par le signe $\ast$ l'une des singularit\'es simples 
d'intersection compl\`ete de courbe espace 
i.e. $\ast = S_{\mu}, U_{\mu}, T_{\mu}, W_{\mu}, 
Z_{\mu}\;\;(\mu \geq 5).$ 
 
Pour chaque singularit\'e $\ast,$ on d\'esigne par $d_j$ la coordonn\'ee 
du point d'intersection du diagramme de Newton du facteur irr\'eductible du 
discriminant $\Delta(s)$ avec l'axe $z_j, j= 1,2.$ 
 
Par r\'ecurrence, on obtient de (4.1.1) une EDF pour $M_k(z_1, z_2).$ 
Quant aux cas $S_{\mu},$ nous renvoyons les lecteurs \`a \cite{AT}. 
 
\begin{lem} 
Soit $\ast$ une des singularit\'es de la liste de Giusti (i.e. SISIC 
courbe espace). Alors pour 
chaque singularit\'e $\ast,$ 
la transform\'ee  de  Mellin  $M_k(z_1,  z_2)$  de  l'int\'egrale 
-fibre 
 $I_k(s_1, s_2)$ 
avec 
 $\nu_1+ \nu_2 + \cdots +\nu_m+1 \leq k \leq 
\nu_1+ \cdots +\nu_{m+1},$ satisfait l'EDF suivante: 
$$M_k(z_1+ d_1, z_2-d_2) 
= V_{m} \prod_{j=k}^{\nu_1+ \cdots +\nu_{m+1}}\frac{1}{L(\ast, k,j;z_1)} 
\times$$ 
$$\times 
\prod_{j= \nu_1+ \cdots +\nu_{m}}^{k-1}\frac{1}{{\tilde L}(\ast, k,j;z_1)} 
\times \prod_{j=1}^{\nu_{m+1}} (z_1 +j ) M_k(z_1, z_2), \leqno (4.1.2)$$ 
o\`u $V_{m} = \prod_{j= \nu_1+ \nu_2 + \cdots +\nu_m +1} 
^{\nu_1+ \nu_2 + \cdots +\nu_m+\nu_{m+1}} v_j,$ produit d'\'el\'ements
de la matrice diagonale $V= diag(v_1, \cdots, v_\mu)$ de la liste de \S 2.
$$ L(\ast, k,j;z_1) = z_1 + \lambda_j +(j-k) + \delta_{k,j}d_1,\;\;k \leq 
j \leq \nu_1+ \cdots +\nu_{m+1},$$ 
$$ {\tilde L}(\ast, k,j;z_1) = z_1 +\lambda_j+ (j-k) + \nu_{m+1}+1, \;\; 
\nu_1+ \cdots +\nu_m+1 \leq j \leq k-1,$$ 
avec $\delta_{k,j},$ delta de Kronecker. Les donn\'ees $\lambda_j, d_1, d_2$ 
sont class\'ees pour chaque  singularit\'e $\ast$ dans la liste du \S 2. 
\label{lem2} 
\end{lem} 
 
Pour r\'esoudre l'\'equation (4.1.2), on \'etablit un lemme. 
\begin{lem} 
L'\'equation aux diff\'erences finies 
$$ M(z + \alpha) = c \cdot \prod_{j=1}^{\nu}\frac{(z+ \beta_j)}{(z+ 
\gamma_j)} M(z),$$ 
admet une solution comme suit: 
$$ M(z ) = c^{\frac{z}{\alpha}} \cdot \prod_{j=1}^{\nu} 
\frac{\Gamma(-\frac{z}{\alpha}-\frac{\gamma_j}{\alpha} ) 
\Gamma(-\frac{\beta_j}{\alpha})}{\Gamma(-\frac{z}{\alpha}- 
\frac{\beta_j}{\alpha}) \Gamma(-\frac{\gamma_j}{\alpha}) } g(z),$$ 
avec une fonction p\'eriodique $ g(z) =g(z + \alpha).$ 
\label{lem3} 
\end{lem} 
 
A l'aide du lemme~\ref{lem3}, on obtient une expression explicite 
de $M_k(z).$ Avant de formuler le th\'eor\`eme, introduisons la notation 
$M_{\phi_k, \gamma_q}(z), 1 \leq k,q \leq \mu$ qui correspond \`a la 
transform\'ee de Mellin 
de l'int\'egrale  $I_{\phi_k,\gamma_q}(s_1,s_2)$ prise le long d'un cycle 
\'evanescent $\gamma_q$ que l'on n'a pas precis\'e dans (1.8.1): 
$$ I_{\phi_k, \gamma_q}(s_1,s_2) 
 =(\frac{1}{2\pi i})^2 \int_{\partial \gamma_q(s)}\frac 
{\phi_k dx_1\wedge dx_2\wedge dx_3 }{(f_1 -s_1)(f_2- s_2)}.$$ 
Et 
$$  M_{\phi_k, \gamma_q} (z_1,z_2) = \int_{\gamma} I_{\phi_k, 
\gamma_q}(s_1,s_2) s_1^{z_1} s_2^{z_2} 
ds_1 ds_2$$ pour un certain 2-cycle r\'eel $\gamma$ qui \'evite les p\^oles de 
$I_{\phi_k, \gamma_q}(s).$ 
\begin{thm} 
Dans la situation d\'ecrite ci-dessus, pour $\nu_1+ \nu_2 
+ \cdots +\nu_m +1 \leq 
k \leq \nu_1+ \cdots +\nu_{m+1},$ on a l'expression suivante: 
\begin{flushleft} 
$$M_{\phi_k, \gamma_q}(z) 
= V_{m}^{\frac{z_1}{d_1}} 
\times  \prod_{j=k}^{\nu_1 + \cdots +\nu_{m+1}} 
\Gamma (- \frac{1}{d_1}L(\ast, k,j;z_1)  ) $$ 
$$\times \prod_{j= \nu_1 + \cdots +\nu_{m}}^{k-1} 
\Gamma (- \frac{1}{d_1}{\tilde L}(\ast, k,j,;z_1)  ) 
\times \prod_{j=1}^{\nu_{m+1}} 
\Gamma ( - \frac{1}{d_1}(z_1 +j ))^{-1} 
\times (\frac{z_2 -\zeta_q}{d_2}+ \frac{z_1}{d_1})^{-1} 
g(z_1) 
, \leqno(4.1.3)$$ 
\end {flushleft} 
avec  $\zeta_q \in {\bf Z}_{\geq 0}, 1 \leq q 
\leq \mu, $ et $\zeta_1 =0. $ 
Ici $g(z_1)$ est une fonction m\'eromorphe
p\'eriodique telle que $g(z_1) = g(z_1 +d_1).$ 
\label{thm43} 
\end{thm} 
{\bf D\'emonstration} 
Il suffit d'appliquer lemme~\ref{lem3} \`a (4.1.2).  
Quant aux  $\zeta_q,$ ils sont d\'etermin\'es par l'\'equation d\'eterminante 
pour les exposants caract\'eristiques \`a $s=0$ du syst\`eme (1.10.1), 
(1.10.2) 
(la m\'ethode de Frobenius) \cite{AK}. Il est facile de voir que la s\'erie 
 des exposants caract\'eristiques  $\zeta_q  $  contient  la s\'erie  des 
entiers  $\zeta_q = 0,1,2,\cdots.$ {\bf C.Q.F.D.} 
 
Nous notons $supp( M_{\phi_k, \gamma_q})(z_1,z_2)$ $\subset$ 
$\frac{1}{d_2} {\bf Q}_{z_1}\times \frac{1}{d_1} {\bf Q}_{z_2} $ 
des points d'intersection  des droites polaires 
$$\frac{1}{d_1}L(\ast, k,j;z_1),\frac{1}{d_1}{\tilde L}(\ast, k,j;z_1), 
 d_1 z_2 + d_2 z_1  \in {\bf Z}_{\geq 0}.$$ 
En liaison avec le Th\'eor\`eme~\ref{thm1}, nous remarquons que 
$$ \frac{d_2}{d_1} = \frac{p_1}{p_2}.$$ 
C'est les points de $supp M_{\phi_k, \gamma_q}(z_1,z_2)$ qui vont 
essentiellement contribuer \`a la transformation inverse de Mellin de 
$M_{\phi_k, \gamma_q}(z_1,z_2)$ qui nous permet de r\'ecup\'erer 
$I_{\phi_k, \gamma_q}(s_1,s_2).$ 
Si on choisit la fonction m\'eromorphe p\'eriodique $g(z_1)$ dans (4.1.3)
en sorte que la transformation inverse de Mellin de 
$M_{\phi_k, \gamma_q}(z_1,z_2)$ ait sens (cf. l'astuce de N\"orlund de 
\S 4.2 ci-dessus), alors on verra facilement la propri\'et\'e suivante de cet 
ensemble. Ici  on fait attention \`a l'in\'egalit\'e $\nu_j \leq d_1 , j = 
1,2,... $ 
 
\begin{cor} Pour un cycle \'evanescent quelconque 
$\gamma,$ 
l'ensemble  $supp (M_{\phi_k, \gamma})(z_1,z_2)$ 
consiste  en les points de la forme 
$$ (  -\lambda_i + a,\frac{d_2}{d_1}(\lambda_i + b) ), 
\;\; 1 \leq i \leq \mu,\;(a,b)\in {\bf Z}^2$$ 
qui sont contenus dans un c\^one $\Gamma_k$ 
$$\Gamma_k = \{ (z_1,z_2) \in{\bf C}^2 ; z_1 \geq -\lambda_k -d_1, 
\frac{z_1}{d_1} + \frac{z_2}{d_2} \geq 0 \}$$ 
de sommet $(-\lambda_k -d_1, \frac{d_2}{d_1}(\lambda_k + d_1) ).$ 
 
\label{cor1} 
\end{cor} 
Si on regarde plus pr\'ecisement  la projection $proj_{z_1} (M_{\phi_1, 
\gamma} 
(z_1,z_2))$ sur l'axe $z_1$ 
de $supp( M_{\phi_k, \gamma})(z_1,z_2)$, les suites suivantes 
s'obtiennent: 
$$ proj_{z_1} (M_{\phi_1, \gamma}(z_1,z_2))
 \subseteq \{ -\lambda_1 -d_1, -\lambda_2 
-1, 
-\lambda_3 -2, \cdots, 1-\nu_1 - \lambda_{\nu_1} \} + d_1{\bf Z}_{\geq 0},$$ 
$$ proj_{z_1} (M_{\phi_2, \gamma}(z_1,z_2))
\subseteq \{ -\lambda_2 -d_1, -\lambda_3 -1, 
-\lambda_4 -2, \cdots, 1-\nu_1 - \lambda_{1} \} + d_1{\bf Z}_{\geq 0}.$$ 
$\vdots$ 
$$ proj_{z_1} (M_{\phi_k, \gamma}(z_1,z_2))\subseteq \{ -\lambda_k -d_1, 
-\lambda_{k+1} -1,-\lambda_{k+2} -2, \cdots, 1-\nu_m - \lambda_{[k-1]} \} + 
d_1{\bf Z}_{\geq 0},$$ 
pour $k$ tel que
 $\nu_1+ \nu_2 + \cdots +\nu_m \leq k \leq 
\nu_1+ \cdots +\nu_{m+1}.$ 
La notation $[k-1]$ est la m\^eme qu'au 
d\'ebut de la section. 
 
\begin{remark} 
{\em  
Il est opportun d'\'evoquer ici le travail \cite{Var} 
de A.N.Varchenko qui a 
d\'efini, dans le cas des singularit\'es isol\'ees d'hypersurface 
associ\'ees au germe $ f: ( {\bf C}^n_x, 0) \rightarrow ( {\bf C}, 0), $ 
 le poids 
d'une forme holomorphe $\omega$ de  la fa\c{c}on suivante: 
$$ \alpha(\omega)= inf \{ \alpha \in {\bf Q}; \frac{1}{t^{\alpha}} 
(\int_{\gamma(t)} \frac{\omega}{df} ) \rightarrow 0 \;\; {\rm lorsque}\; 
t \rightarrow 0 \}$$ 
$$ = min\{ \alpha \in {\bf Q}; \mbox{il existe } \;0\leq k 
\leq n-1 \;\mbox{ 
tel que}\; A_{k,\alpha} \not= 0$$ 
$$\mbox {du d\'eveloppement asymptotique}\; 
 \int_{\gamma(t)} \frac{\omega}{df}  \sim \sum_{\lambda \in \Lambda} 
\sum_{\alpha \in L(\lambda)} \sum _{0 \leq k \leq n-1}   A_{k,\alpha} 
t^{\alpha} (\log\; t)^{k} \},$$ 
o\`u $\gamma(t)$ est un cycle \'evanescent lorsque 
$t \rightarrow 0,$  $L(\lambda)= 
\{\beta >-1; exp( -2\pi i \beta) = \lambda  \},$ 
$\Lambda=$ $\{$  les valeurs propres  $\lambda\;$ de la monodromie de 
Picard-Lefschetz de la singularit\'e   $ f(x)=0  \;\;\}.$ 
 
Il est facile de voir que dans notre situation le bord du c\^one 
$ \Gamma_k$ introduit dans le Corollaire ~\ref{cor1} correspond \`a 
$\alpha(\omega)$ de Varchenko. En particulier, si la valeur critique de
$f$ est constitu\'ee d'un cusp (et d'une droite en position g\'en\'erique)
notre $\lambda_i$ donne la monodromie de Picard-Lefschetz de la singularit\'e.}
\label{remark411}
\end{remark}

A juste titre, on pourrait donner une autre d\'efinition plus g\'en\'erale
des spectres 
du syst\`eme de Gauss-Manin au lieu de la D\'efinition ~\ref{dfn31}. 
 
\begin{dfn}
Les spectres du syt\`eme de Gauss-Manin associ\'e aux SIIC courbe-espace
quasihomog\`ene  consistent en les  donn\'ees suivantes: 
la partie de l'ensemble des droites 
contenue dans le bord d'un c\^one $ \bigcup_\gamma$
$supp(M_{\phi_k, \gamma})(z_1,z_2),$ $1 \leq k \leq \mu,$
o\`u $\gamma $ parcourt tout les cycles \'evanescents de la singularit\'e. 
\label{dfn411} 
\end{dfn} 
 
D'apr\`es cette nouvelle d\'efinition, les spectres du syst\`eme (3.1.1)
sont donn\'es par le bord d'ensembles
 $(\{ z_1\geq -\lambda_k -d_1\} \cap \{ d_1 z_2 + d_2 z_1 \geq 0\}),$
$1 \leq k \leq \mu.$
 
\begin{remark} 
{\em 
Il faut remarquer que les r\'eseaux de l'ensemble 
$supp(M_{\phi_k, \gamma})(z_1,z_2)$ donnent naissance
\`a la bonne $\kappa$-filtration ($\kappa$ =2) introduite par C.Sabbah 
\cite{Sab1}, \cite{Sab2}. 
 
Donc il suffit d'appliquer sa Proposition 1.2. de \cite{Sab2} \`a cette 
bonne filtration, 
pour voir l'existence de l'ensemble ${\cal L}_1$ de formes lin\'eaires  et 
d'un polyn\^ome \`a une variable $b_{L,1}$ tels que 
$$ [\prod_{L \in {\cal L}_1}b_{L,1}(L(s)) ]f_1^{s_1}f_2^{s_2} 
= P_1(x_1,x_2,x_3;\frac{\partial}{\partial x_1},\frac{\partial}{\partial x_2}, 
\frac{\partial}{\partial x_3},s_1, s_2 )  f_1^{s_1+1}f_2^{s_2},$$ 
pour les polyn\^omes $f_1, f_2$ trait\'es dans cet article. 
Ici nous avons comme l'ensemble ${\cal L}_1,$ 
${\cal L}_1 =\{ s_1 +\lambda_k + d_1 -\alpha, d_1 s_2 + d_2 s_1 - 
\gamma \},$ avec certains
$\alpha, \gamma \in {\bf Z_{\geq 0}},$ 
et $\lambda_1, \cdots, \lambda_\mu$ les spectres du syst\`eme de 
Gauss-Manin. Il est naturel de consid\'erer le polyn\^ome 
$\prod_{L \in {\cal L}_1}b_{L,1}(L(s))$ ci-dessus
comme la  $b-$fonction en 
2-variables.
 
Quant \`a l'EDF correspondant aux d\'ecalages en $s_2,$ elle s'obtient 
d'une fa\c{c}on analogue en partant du syst\`eme (3.1.2). 
 
$$ [\prod_{L \in {\cal L}_2} b_{L,2}(L(s)) ]f_1^{s_1}f_2^{s_2} 
= P_2(x_1,x_2,x_3;\frac{\partial}{\partial x_1},\frac{\partial}{\partial x_2}, 
\frac{\partial}{\partial x_3},s_1, s_2 )  f_1^{s_1}f_2^{s_2 +1},$$ 
o\`u 
${\cal L}_2 =$ $\{ s_2 -\frac{d_2}{d_1}\lambda_k -\beta,$ 
$d_1 s_2 + d_2 s_1 -\gamma \},\;$ avec certains
 $\beta, \gamma \in {\bf Z_{\geq 0}}.$ Les notations sont les m\^emes que 
celles de  ${\cal L}_1.$ 
} 
\label{remark412} 
\end{remark} 
 
{\bf 4.2. L'int\'egrale- fibre en tant qu'une fonction hyperg\'eom\'etrique 
g\'en\'eralis\'ee } 
 
D'ailleurs il serait utile de voir le r\'esultat du Th\'eor\`eme 
~\ref{thm43} en liaison avec la notion  des fonctions 
hyperg\'eom\'etriques g\'en\'eralis\'ees (FHG) au sens de 
Mellin-Barnes-Pincherle 
\cite{AK}, \cite{Nor}. 
Par cette formulation la FHG de Gauss s'exprime par l'int\'egrale, 
$$\frac {1}{2\pi i } \int_{z_0 - i\infty}^{  z_0 + i\infty }(-s)^z 
\frac {\Gamma ( z+ \alpha)\Gamma ( z+ \beta )\Gamma ( -z)}{\Gamma ( z+ 
\gamma)} 
dz , \;\; - \Re  \alpha, - \Re  \beta < z_0.$$ 
 
Pour assurer la convergence de la transform\'ee de Mellin inverse de 
$M_{\phi_k, \gamma}(z)$ de (4.1.3): 
$$ \int_{\Pi}s_1^{-z_1-1} s_2^{-z_2-1}M_{\phi_k, \gamma}(z) 
dz_1 dz_2 , \leqno(4.2.1)$$ 
on v\'erifie que l'EDF (4.1.2) admet la solution $M_{\phi_k}(z)$ 
telle que pour un certain $\epsilon >0,$ 
$$\mid M_{\phi_k}(z)\mid < C_k exp(-\epsilon \mid Im\; z \mid )\;\;{\rm 
lorsque}\;  Im\; z \rightarrow \infty, 
\mbox{ dans un secteur d'ouverture }\;<2 
\pi. $$ 
 
Pour voir l'existence d'une solution de l'EDF avec d\'ecroissance 
exponentielle, on recourt \`a une astuce de N\"orlund \cite{Nor}. 
Sa technique consiste en un choix du facteur $g(z_1)$ de (4.1.3) 
qui doit \^etre 
une fonction m\'eromorphe de p\'eriode $d_1.$ 
Si on note $z= -\frac{z_1}{d_1},$ notre analyse de (4.2.1) est reduite \`a 
l'\'etude 
de l'int\'egrale 
$$ \int_{z_0 - i\infty}^{  z_0 + i\infty } s^z g(z) 
\prod_{j=1}^{\nu}\frac {\Gamma ( z+ \alpha_j)}{\Gamma ( z+ \rho_j)} 
dz . \leqno(4.2.2)$$ 
 
\begin{lem} 
Si on choit une des fonctions suivantes $g^{+}(z)$ (resp. $g^{-}(z)$) en tant 
que $g(z),$ 
alors l'int\'egrand de (4.2.2) est de d\'ecroissance exponentielle lorsque 
$ Im \;z $ tend vers $\infty $ dans le secteur $0 \leq arg \;z < 2\pi,$ 
(resp. $-\pi \leq arg \;z < \pi .$) 
 
$$ g^{\pm}(z)= 1+ e ^{\pm 2 \pi i \beta_{\nu}}\prod_{j=1}^{\nu} 
\frac {sin 2 \pi ( z+ \alpha_j)}{sin 2\pi ( z+ \rho_j)}, $$ 
avec $ \beta_{\nu} =-1 + \sum_{j=1}^{\nu}(\rho_j  - \alpha_j)  $ 
\label{lem45} 
\end{lem} 
 
{\bf D\'emonstration} 
 
Il suffit de se rappeler 
$$\prod_{j=1}^{\nu}\frac {\Gamma ( x+ iy + \alpha_j)}{\Gamma ( x+iy + \rho_j)} 
\rightarrow const. \mid y \mid^{-(\beta_{\nu}+1)}$$ 
lorsque  $y \rightarrow \pm \infty.$ Ici, on se sert de la formule 
de Binet: 
$$ log \Gamma(z+a) = \log \Gamma(z) + a \log z - \frac{a- a^2}{2z} + {\cal O} 
( \mid z \mid^{-2}   ) 
$$ 
si $\mid z \mid >> 1,$ (Whittaker-Watson, Chapter XII, Example 44). 
Le facteur $\mid s^{-(x+iy)} \mid = r^{-x}e^{\theta y},$ 
pour $s =  re^{i \theta} $ donne la contribution exponentiellement 
d\'ecroissante dans chaque cas. 
{\bf C.Q.F.D.} 
 
Ainsi on a d\'emontr\'e la convergence de l'expression (4.2.1) 
pour un certain $\Pi$ qui est obtenu comme un produit du chemin 
d'int\'egration dans ${\bf C}_{z_1}$ et celui dans ${\bf C}_{z_2}.$ 
\begin{prop} 
L'int\'egrale- fibre 
$I_{\phi_k, \gamma_q}(z)$  des SISIC  prise le long d'un cycle 
\'evanescent $\gamma_q$ 
est une FHG au sens 
de Mellin-Barnes-Pincherle d\'efinie par l'expression suivante 
$$ \int_{\Pi} s_1^{-z_1-1}s_2^{-z_2-1}M_{\phi_k, \gamma_q}(z) 
dz_1 dz_2,$$ 
pour le $ M_{\phi_k, \gamma_q}(z)$ qui appara\^{\i}t en 
(4.1.3) avec $g(z_1)= g^+(-\frac{z_1}{d_1})$ 
introduit dans Lemme ~\ref{lem45} 
 
\label{prop12} 
\end{prop} 
\vspace{2pc} 
{ 
\center{\section{ 
Des cas unimodaux et des autres cas accessibles }} 
} 
  Malgr\'e le caract\`ere restrictif de calculs faits dans les sections 
pr\'ecedentes, nos d\'emarches s'appliquent aux autres cas qui contiennent 
des s\'eries infinies. 
 
{\bf 5.1. Le cas non-r\'esonant.} 
 
On regarde l'application quasihomog\`ene suivante, 
$$ \leqno(5.1.1)$$  
$$ 
\left\{ 
\begin{array}{ccccc} 
f_1(x_1,x_2,x_3) & = & x_1^{q_1} + x_2^{q_2} +x_3^{q_3} &=  & 0\\ 
f_2(x_1,x_2,x_3) & = & x_1x_2 & = &0 \\ 
\end{array} 
\right. 
$$ 
o\`u on impose la condition que les poids $w_i$ et $w_j$ ( 
$ 1 \leq i,j \leq 3 $) soient premiers entre eux. Nous l'appelons le 
cas non-r\'esonant. 
Selon la notation du \S 2, $p_1 = q_1w_1= 
q_2 w_2 = q_3w_3, p_2 = w_1 + w_2 .$ 
Pour voir l'analogie du cas (5.1.1) avec les cas de singularit\'es isol\'ees 
simples de courbe espace, on \'etablit l'enonc\'e suivant: 
 
\begin{lem} 
   Pour $(f_1,f_2)$ de (5.1.1) sous l'hypoth\`ese sur les poids $w_1,w_2,w_3$ 
comme ci-dessus, les matrices $P^{(1)}, P^{(2)}$ d\'efinies dans 
la Proposition ~\ref{prop1} 
sont semblables \`a des matrices diagonales. 
Plus pr\'ecisement, pour chaque $1 \leq j \leq \mu,$ il 
existe au plus des uniques $j_1$ et $j_2$ tels que  
$$ \tilde{\omega_j}\wedge df_2 =  P^{(1)}_{jj_2} \phi_{j_2} dx_1\wedge dx_2\wedge dx_3 \leqno(5.1.2)$$ 
$$ \tilde{\omega_j}\wedge df_1 =  P^{(2)}_{jj_1} \phi_{j_1} dx_1\wedge dx_2\wedge dx_3 \leqno(5.1.3)$$ 
avec $P^{(1)}_{jj_2}\not = 0, P^{(2)}_{jj_1} \not =0$ et $j_1
\not = j_2.$ En revanche pour chaque indice $j$, on trouve au plus 
des entiers uniques $\tilde{j_1}, \tilde{j_2}$ tels que
$$
\tilde \omega_{\tilde j_2} \wedge df_1  =  P^{(2)}_{\tilde{j_2}j}\phi_j(x)dx_1\wedge dx_2\wedge dx_3 
\leqno(5.1.4)$$
$$\tilde \omega_{\tilde j_1} \wedge df_2  =  P^{(1)}_{\tilde{j_1}j}\phi_j(x)dx_1\wedge dx_2\wedge dx_3 \leqno(5.1.5)$$
avec $P^{(2)}_{\tilde{j_2} j} \ne 0, P^{(1)}_{\tilde{j_1} j} \ne 0$. 
\label{lem8} 
\end{lem} 
 
{\bf D\'emonstration} 
 
  Nous faisons la comparaison entre les poids  des termes. 
On d\'emontre le cas (5.1.2). Le cas (5.1.3) se d\'emontre 
d'une mani\`ere similaire.
 
S'il existe un 
autre terme \`a part de $ \phi_{j_2}$ (disons  $\phi_{j_2'}$ ) qui participe 
\`a la d\'ecomposition (5.1.2), leurs poids  satisfont les 
relations suivantes: 
 $$ \ell_j + p_2 = w( P^{(1)}_{jj_2}) + w (\phi_{j_2})+ w_1 +w_2 +w_3$$ 
$$ = w( P^{(1)}_{jj_2'}) + w (\phi_{j_2'})+ w_1 +w_2 +w_3.    
\leqno(5.1.6)$$ 
Remarquons que $w( P^{(1)}_{jj_2}),  
w( P^{(1)}_{jj_2'}) \in p_1 {\bf Z}_{\geq 0} 
+ p_2 {\bf Z}_{\geq 0}.$ 
   Donc la diff\'erence des poids  de deux formes $w (\phi_{j_2})$ 
et $ w (\phi_{j_2'})$ doit appartenir au r\'eseau des entiers 
$p_1 {\bf Z} + p_2 {\bf Z}.$ 
 
 Maintenant on se souvient de la formule de la s\'erie de Poincar\'e de
l'espace  $\tilde \Phi,$ (3.2.6). Pour le cas (5.1.1), elle devient, 
\begin{flushleft} 
$ P_{\tilde \Phi}(t)=$ 
\end{flushleft} 
$$= t^{(q_3-1)w_3}+ (1+ t^{w_3}+ \cdots +t^{(q_3-2)w_3} ) 
(1+ t^{w_2}+ \cdots +t^{q_2w_2} + t^{w_1}+ \cdots +t^{(q_1-1)w_1}). 
 \leqno(5.1.7) 
$$ 
De la formule (5.1.7), il est facile de voir que $w (\phi_{j_2})- 
 w (\phi_{j_2'}) \not \in  p_1 {\bf Z}_{\geq 0} 
+ p_2 {\bf Z}_{> 0} = q_i w_i {\bf Z}_{\geq 0} 
+(w_1 + w_2) {\bf Z}_{> 0}. $ 
Cela veut dire que $\frac {P^{(1)}_{jj_2}}{P^{(1)}_{jj_2'}} \not \in s_2 {\bf C}[s]$ 
si $P^{(1)}_{jj_2'} \not = 0.$ Autrement dit, si on trouve deux  indices $j_2, 
j_2' $ pour lesquels (5.1.6) sont v\'erifi\'e, alors ils sont pour ces indices:
$$ w (\phi_{j_2}) = kw_3 + p_1, w (\phi_{j_2'})= kw_3, 0 \leq k \leq q_3-1,$$ 
d'apr\`es (5.1.7).  Cette relation entra\^{\i}ne que $P^{(1)}_{jj_2}= cf_1 P^{(1)}_{jj_2'}$ 
pour certain constant $c:$ 
$$\tilde{\omega_j}\wedge df_2 = (P^{(1)}_{jj_2'} \phi_{j_2'} + 
P^{(1)}_{jj_2} \phi_{j_2}) 
dx_1\wedge dx_2\wedge dx_3 $$ 
$$ = P^{(1)}_{jj_2'}( \phi_{j_2'} + cf_1 \phi_{j_2})dx_1\wedge dx_2\wedge dx_3. \leqno(5.1.8)$$ 
Il faut remarquer ici que le terme de gauche est un mon\^ome pour 
les singularit\'es d\'efinies par (5.1.1), par contre le 
terme de droite doit \^etre essentiellement polynomial sinon $c=0.$ 
Donc il faut 
qu'un seul terme $P^{(1)}_{jj_2'} \phi_{j_2'}$ prenne part \`a 
la d\'ecomposition 
(5.1.8).
  
 Pour voir que les indices $j_1, j_2$
de (5.1.2), (5.1.3) sont diff\'erents, il faut calculer l'espace
$F:$
\begin{flushleft}
$ P_{F}(t)=$  
\end{flushleft}
$$(1+ t^{w_3}+ \cdots +t^{(q_3-1)w_3} )t^{w_1 +w_2} +
t^{w_3}(1+ t^{w_3}+ \cdots +t^{(q_3-2)w_3}) (t^{w_2}+ 
\cdots +t^{q_2w_2} + t^{w_1}+ \cdots +t^{(q_1-1)w_1}). 
$$ 
Des formes qui pourraient produire la situation avec $j_1= j_2$
dans (5.1.2),(5.1.3) sont celles dont les  poids appartiennent
au r\'eseau $w_1+w_2+w_3 + w_1 {\bf Z}_{\geq0} + 
w_2{\bf Z}_{\geq0}  + w_3 {\bf Z}_{\geq0}.$ Par le 
calcul direct des formes de $F,$ les uniques formes qui satisfont cette condition sont $x_2x_3^i dx_3dx_1 \in F\;$ 
$ (0 \leq i \leq q_3-2).$ Cela ach\`eve la d\'emonstration
de (5.1.2) et d'une fa\c{c}on analogue (5.1.3). D'apr\`es
une comparaison des s\'eries de Poincar\'es $ P_{F}(t)$ et
$P_{\tilde \Phi}(t)$ on peut conclure (5.1.4),(5.1.5) en 
tenant compte de (5.1.2),(5.1.3). 
{\bf C.Q.F.D. } 
 
On remarque ici que la condition  de 
non-r\'esonance a \'et\'e essentiellement utilis\'ee pour que le terme de 
gauche
de (5.1.8) soit monomial.  
 
{\bf 5. 2. La liste de Wall et d'Aleksandrov } 
 
Pour les singularit\'es d'intersection compl\`ete pas n\'ecessairement 
simples, les r\'esultats analogues aux Th\'eor\`emes ~\ref{thm1}, ~\ref{thm43} 
et ~\ref{cor1} sont valables. 
Notamment la s\'erie des singularit\'es unimodales 
de la liste de Wall \cite{Wa}: 
 
$P_{k, \ell}, k\geq \ell \geq 2, \mu= k + \ell + 1,$ 
\begin{flushleft} 
 
\[ 
\left\{ 
\begin{array}{ccccc} 
f_1(x_1,x_2,x_3) & = & x_1^k + x_2^{\ell} + x_3^2 & =  & 0\\ 
f_2(x_1,x_2,x_3) & = & x_1x_2  & = &0 \\ 
\end{array} 
\right. 
\] 
 
\end{flushleft}

$G_{2m+6}, m\geq 3$ 
\begin{flushleft}  
\[ 
\left\{ 
\begin{array}{ccccc} 
f_1(x_1,x_2,x_3) & = & x_1^2+ x_2x_3^m & =  & 0\\ 
f_2(x_1,x_2,x_3) & = & x_2^2 + x_3^3 & = &0 \\ 
\end{array} 
\right. 
\] 
\end{flushleft}

$G_{2m+3},m\geq 3$ 
\begin{flushleft} 
 
\[ 
\left\{ 
\begin{array}{ccccc} 
f_1(x_1,x_2,x_3) & = & x_1^2+ x_3^m & =  & 0\\ 
f_2(x_1,x_2,x_3) & = & x_2^2 + x_3^3 & = &0 \\ 
\end{array} 
\right. 
\] 
 
\end{flushleft} 
\begin{flushleft}
$FZ_{6m+6}$

\[
\left\{
\begin{array}{ccccc}
f_1 & = & x_1x_3 + x^3_3 + x^{3m+1}_2 \\ 
f_2 & = & x_1x_2 
\end{array}
\right.
\]
\end{flushleft} 
\begin{flushleft}
$FZ_{6m+8}$
\[
\left\{
\begin{array} {ccccc}
f_1 & = & x_1x_3+x^3_3 + x^{3m+2}_2 \\ 
f_2 & = & x_1x_2\\
\end{array}
\right.
\]
\end{flushleft} 

 On remarque  la structure du module $F$ analogue \`a celle 
des cas 
$S_\mu, T_\mu$ ( pour $P_{k,\ell}$) et des cas $Z_{10}$ ($G_{2m+3}$ et 
$G_{2m+6}$). 
C'est-\`a- dire: 
\begin{flushleft} 
$P_{k,\ell}$ 
$$F  = \{ dx_1dx_2 , x_2dx_3dx_1,$$
$$
\overbrace{ dx_2dx_3, x_2dx_2dx_3,\cdots, 
x_2^{\ell -2}dx_2dx_3}^{\ell-1},  
\overbrace{dx_3dx_1, x_1dx_3dx_1, \cdots,
x_1^{k-1}dx_3dx_1}^{k} \}.$$ 
$$  w_1 = 2\ell, w_2 = 2k, w_3 =k\ell, p_1 = 2k\ell, p_2= 2 (k+ \ell)w. $$ 
$$\tilde \Phi = \{ 1, x_1, \cdots, x_1^{k-1}, x_2, \cdots, x_2^{\ell}, x_3 \}.$$ 
$$w(\phi_j)=\{ 0, 2\ell, \cdots ,2(k-1)\ell, 2k, \cdots ,2k\ell, k\ell \}.$$

\end{flushleft} 
 
\begin{flushleft} 
$G_{2m+6}$ 
$$F  = \{dx_2dx_3,x_3dx_2dx_3, dx_3dx_1, dx_1dx_2,x_3dx_3dx_1,$$ 
$$\overbrace{ 3x_2dx_3dx_1+2x_3dx_1dx_2, x_3(3x_2dx_3dx_1+2x_3dx_1dx_2),
\cdots,
x_3^{m-2}(3x_2dx_3dx_1+2x_3dx_1dx_2)}^{m-1},  
$$ $$\overbrace{
x_2( mx_1dx_2dx_3+2x_3dx_1dx_2),x_2x_3( mx_1dx_2dx_3+2x_3dx_1dx_2),\cdots, 
x_2x_3^{m-2}(mx_1dx_2dx_3+2x_3dx_1dx_2)}^{m-1},$$
$$mx_1dx_2dx_3+2x_3dx_1dx_2, x_3( mx_1dx_2dx_3+2x_3dx_1dx_2),mx_3^2dx_3dx_1+x_2dx_1dx_2  
\}.$$ 
$$  w_1 = 2m+3, w_2 = 6, w_3 = 4, p_1 = 2(2m+3), p_2= 12. $$ 
$$\tilde \Phi = \{ 1, x_3,x_3^2, \cdots, x_3^{m+1}, x_2,x_2x_3 \cdots, x_2x_3^{m+1},x_1,x_1x_3 \}.$$ 
$$w(\phi_j)=\{ 4i, 6+ 4i \; ( 0 \leq i \leq m+1), 2m+3, 2m+7 \}.$$ 

\end{flushleft}

\begin{flushleft} 
$G_{2m+3}$ 
$$F = \{dx_2dx_3,x_3dx_2dx_3, dx_1dx_2,x_3dx_1dx_2, x_3^2dx_1dx_2, 
\overbrace{ dx_3dx_1, x_3dx_3dx_1,\cdots, 
x_3^{m-2}dx_3dx_1}^{m-1},$$  
$$\overbrace{ mx_1dx_2dx_3+2x_3dx_1dx_2, 
x_3( mx_1dx_2dx_3+2x_3dx_1dx_2),\cdots, 
x_3^{m-2}(mx_1dx_2dx_3+2x_3dx_1dx_2)}^{m-1}
\}.$$
$$  w_1 = m, w_2 = 3, w_3 = 2, p_1 = 2m, p_2= 6. $$ 
$$\tilde \Phi = \{ 1, x_3,x_3^2, \cdots, x_3^{m+1}, x_2,x_2x_3 \cdots, x_2x_3^{m-2},x_1,x_1x_3 \}.$$ 
$$w(\phi_j)=\{ 2i\; ( 0 \leq i \leq m+1), 3+ 2j\; (0 \leq j \leq m-2),m, m+2 
\}.$$ 
\end{flushleft} 
 
\begin{flushleft}
$FZ_{6m+6}$

$$\begin{array}{c}
F= \{ dx_3dx_1, x_3dx_2dx_3, x^i_2(3x_3dx_2dx_3 - dx_1dx_2), 
0 \le i \le 3m-1, x_1(dx_1dx_2 - 3x_3dx_2dx_3),\\ x_1dx_1dx_3, 
x^j_2dx_2dx_3, 0 \le j \le 3m, x_2dx_3dx_1\} \\
w_1 = 2(3m+1), w_2 = 3, w_3 = 3m+1, p_1 = 6m+5, p_2 = 3(3m+1) \\
\tilde \Phi = \{ x_1, x^2_3, x^2_1, x^i_2, 0 \le i \le 3m+1, x_2^kx_3, 0 \le k \le 3m\} \\
w(\phi_1) = \{ 2(3m+1), 4(3m+1), 3i, 0 \le i \le 3m+1, 3k+3m+1,
 0 \le k  \le 3m\}.
\end{array}$$
\end{flushleft}

\begin{flushleft}
$FZ_{6m+8}$

$$\begin{array}{c}
F =\{dx_3dx_1, x_3dx_2dx_3, x^i_2(3x_3dx_2dx_3 - dx_1dx_2), 
0 \le i \le 3m, x_1(dx_1,dx_2 - 3x_3dx_2dx_3), \\ 
x_1dx_1dx_3, x^k_2dx_2dx_3, 0 \le k \le 3m+1, x_2dx_3dx_1\} \\
w_1 = 2(3m+2), w_2=3, w_3 = 3m+2, p_1 = 6m+7, p_2 = 3(3m+2), \\
\tilde \Phi = \{x_1, x^2_3, x^2_1, x^i_2, 0 \le i \le 3m+2, 
x^k_2x_3, 0 \le k \le 3m+1\} \\
w(\phi_j) = \{2(3m+2), 2(3m+2), 4(3m+2), 2i, 0 \le i \le 3m+2, 
2k+3m+2, 0 \le k \le 3m+1\}. 
\end{array}$$
\end{flushleft}

On trouve chez A.G.Aleksandrov\cite{Al0} une s\'erie des singularit\'es 
$U_{2m+1}$ ( dans \cite{Wa} on trouve $J_{6m+7},$ $J_{6m+9}$) 
dont le syst\`eme de Gauss-Manin  peut \^etre calcul\'e d'une fa\c{c}on 
analogue au cas $U_7, U_9:$ 
 \begin{flushleft} 
$U_{2m+1}, m\geq 5$ 
 
\[ 
\left\{ 
\begin{array}{ccccc} 
f_1(x_1,x_2,x_3) & = & x_1x_2 + x_3^m & =  & 0\\ 
f_2(x_1,x_2,x_3) & = & x_1^2 + x_2x_3  & = &0 \\ 
\end{array} 
\right. 
\] 
$$F = \{dx_2dx_3,dx_1dx_2, x_3^{m-1}dx_3dx_1+x_1dx_1dx_2,
\overbrace{ dx_3dx_1, x_3dx_3dx_1,\cdots, 
x_3^{m-2}dx_3dx_1}^{m-1},$$  
$$\overbrace{ x_3dx_2dx_3+x_1dx_3dx_1, 
x_3( x_3dx_2dx_3+x_1dx_3dx_1),\cdots, 
x_3^{m-2}( x_3dx_2dx_3+x_1dx_3dx_1)}^{m-1}
\}.$$ 
$$  w_1 = m+1, w_2 = 2m-1, w_3 = 3, p_1 = 3m, p_2= 2m+2. $$ 
$$\tilde \Phi = \{ 1, x_3,x_3^2, \cdots, x_3^{m-1}, x_1,x_1x_3 \cdots, x_1x_3^{m-1},x_2 \}.$$ 
$$w(\phi_j)=\{ 3i, m+1 +3i\; (0 \leq i \leq m-1) , 2m-1 \}.$$ 

\end{flushleft}

Par un calcul semblable \`a celui des \S 2, \S 3, on obtient 
un analogue des Th\'eor\`emes ~\ref{thm1} et ~\ref{cor1}. On reprend la 
notation,  
$$ I_{k, \gamma}(s) 
 = (\frac{1}{2\pi i})^2\int_{\partial \gamma(s)}\frac 
{\phi_k dx_1\wedge dx_2\wedge dx_3 }{(f_1 -s_1)(f_2- s_2)},$$ 
$$  M_{k, \gamma} (z) = \int_{\Pi} I_{k, \gamma} 
(s_1,s_2) s_1^{z_1} s_2^{z_2} 
ds_1 ds_2,$$ pour un certain $\Pi$ qui \'evite les p\^oles de 
$I_{k, \gamma}(s),$ 
et 
${\bf I}_{\Phi}  = (I_{1, \gamma}(s_1,s_2), \cdots, I_{\mu, \gamma}(s_1,s_2)).$ 
En somme, le calcul des cas trait\'es dans 5.1 et 5.2 donnent le r\'esultat
suivant. 
\begin{thm}

{\bf 1.} 
Le syst\`eme de Gauss-Manin pour ${\bf I}_{\Phi}$ associ\'e 
aux singularit\'es isol\'ees d'intersection 
compl\`ete de courbe espace des cas non-r\'esonants (5.1.1), 
$P_{k,l}, G_{2m+6}, G_{2m+3},$ $FW_{13},$ $FW_{19}, K_{13}$ $FZ_{6m+6},$ 
$FZ_{6m+8},$ $HD_{13}, HD_{14},$ 
$K_{14}$ de la liste 
de Wall et $  U_{2m+1} $ de la liste d'Aleksandrov  s'\'ecrit 
sous la forme suivante: 
$$ P'^{(1)}(s_1)( s_1 id_{\mu}\frac{\partial}{\partial 
s_1}   - \frac{1}{p_1}L_{\Phi}) {\bf I}_{\Phi} = \frac{ p_2}{p_1} s_2 
P'^{(2)}(s_2)\frac{\partial}{\partial s_1}{\bf I}_{\Phi}, \leqno(5.2.1)$$ 
o\`u $P'^{(1)}(s_1):$ une matrice diagonale, et $P'^{(2)}(s_2):$ une matrice 
semblable \`a une matrice diagonale. 
$$ \frac{1}{p_1}L_{\Phi} = diag\{\lambda_1, \cdots, \lambda_{\mu}   \}, 
\lambda_j  = \frac {w(\phi_j) +w_1+w_2+w_3 - p_1 -p_2}{p_1}.$$ 
D'une fa\c{c}on analogue: 
$$P'^{(2)}(s_2) ( s_2 id_{\mu}\frac{\partial}{\partial 
s_2}   - \frac{1}{p_2}L_{\Phi}) {\bf I}_{\Phi} = \frac{ p_1}{p_2} s_1 
P'^{(1)}(s_1)\frac{\partial}{\partial s_2}{\bf I}_{\Phi}, \leqno(5.2.2)$$ 
avec les spectres $\{\lambda_1, \cdots, 
\lambda_{\mu}  \}$ poss\'edant la propri\'et\'e de sym\'etrie. 

{\bf 2.} Soit $\ast$
une des singularit\'es suivantes: cas non-r\'esonants (5.1.1),
$P_{k,l},$ $G_{2m+6}, G_{2m+3},$ $U_{2m+1},$ $FW_{13},$ $FW_{19},$ $K_{13}.$ 
Alors on a la formule suivante de la transform\'ee de Mellin de l'int\'egrale- fibre
$M_{k, \gamma_q(z)}$
associ\'ee \`a 
$\ast$.
Pour chaque 
$\phi_k \in {\tilde \Phi}$,
on trouve des ensembles d'entiers
$I^{(k)}_1 \ne \emptyset$
et
$I^{(k)}_2 \ne \emptyset$
tels que
$I^{(k)}_1 \cap I^{(k)}_2 = \emptyset, I^{(k)}_1 \cup I^{(k)}_2 
\subset [1,\mu], k \in I^{(k)}_1$
et
$$\begin{array}{c}
M_{k,\gamma_q}(z) = V^{\frac{z_1} {d_1(k)}}_k \prod_{i\in I^{(k)}_1} \Gamma (- \frac 1 {d_1(k)} L(\ast, k ,i; z_1)) \times \\
\times \prod_{j \in I^{(k)}_2} \Gamma (-\frac 1 {d_1(k)} \tilde L (\ast,k,j;z_1)) \prod^{\mid I^{(k)}_1\mid + \mid I^{(k)}_2 \mid}_{t=1} \Gamma (-\frac 1 {d_1(k)} (z_1+t))^{-1} \times \\
\times (p_2(z_2-\zeta_q) + p_1z_1)^{-1} g(z_1) 
\end{array}$$
o\`u on utilise les notations ci-dessous:
$$\begin{array}{c}
L(\ast, k,i;z_1) = z_1 + \lambda_i + (i-k) + \delta_{k,i}d_1(k) \\
\tilde L(\ast , k,j;z_1) = z_1 + \lambda_j + (j-k) + \mid I^{(k)}_1 \mid + \mid I^{(k)}_2 \mid + 1, 
\end{array}$$ 
o\`u $\delta_{k,j}$
delta de Kronecker,
$d_1(k), \zeta_q \in {\mathbf Z}_{\geq 0}, V_k \in {\mathbf Q},$ $g(z_1)$
une fonction m\'eromorphe p\'eriodique telle que
$g(z_1 + d_1(k)) = g(z_1)$.
Les rationnels (spectres)
$\lambda_j$
sont d\'efinis comme
$\lambda_j = \frac {w(\phi_j)+ \mid w\mid - \mid p \mid} {p_1}$.

{\bf 3.} Soit
$\ast$
une des singularit\'es suivantes de la liste de Wall: 
$FZ_{6m+6}, FZ_{6m+8},$ $HD_{13}, HD_{14},$ $K_{14}$.
Alors, pour chaque
$\phi_k \in \tilde \Phi$
on trouve des ensembles d'entiers
$I^{(k)}_1 \ne \emptyset, I^{(k)}_2 \ne \emptyset, I^{(k)}_3$
et un entier
$\tilde k \in [1,\mu]$
tels que
$I^{(k)}_1 \cap I^{(k)}_2 = \emptyset, I^{(k)}_2 \cap I^{(k)}_3 = \emptyset$,
$I^{(k)}_1 \cap I^{(k)}_3 = \{\tilde k\}, 
I^{(k)}_1 \cup I^{(k)}_2 \cup I^{(k)}_3 \subset [1,\mu]$
avec lesquels la transform\'ee de Mellin
$M_{k,\gamma_q}(z)$
s'\'ecrit sous la forme suivante:

$$\begin{array}{c}
M_{k,\gamma_q}(z) = V^{\frac{z_1+\tau_k}{d_1(k)}} \prod_{i\in I^{(k)}_1} \Gamma (- \frac 1 {d_1(k)} L(\ast, \tilde k, i; z_1+ \tau_k)) \times \\
\times \prod_{j \in I^{(k)}_2} 
\Gamma (- \frac 1 {d_1(k)} \tilde L (\ast, \tilde k, j; z_1 + \tau_k)) \prod^{\mid I^{(k)}_1 \mid + \mid I^{(k)}_2 \mid}_{t=1} \Gamma(-\frac 1 {d_1(k)} (z_1+t+\tau_k)) \times \\
\times \prod_{r \in I^{(k)}_3 \setminus \{\tilde k\}} 
(\frac{z_1+\tilde{\alpha}_r} {z_1 + \tilde{\beta}_r}) 
(p_1(z_1+\tau_k) + p_2(z_2 - \zeta_q + \sigma_k))^{-1} g(z_1) 
\end{array}$$
o\`u
$\sigma_k, \tau_k \in \mathbf Z$
et soit
$\tilde{\alpha}_r \in {\mathbf Z}, 
\tilde{\beta}_r \in \{\lambda_1, \ldots, \lambda_{\mu} \} + {\mathbf Z}$
soit
$\tilde{\alpha}_r \in \{\lambda_1,\ldots, \lambda_{\mu} \} + {\mathbf Z}, 
\tilde{\beta}_r \in \mathbf Z$.
Ici on utilise la notation
$\{\lambda_1,\ldots, \lambda_{\mu}\} + \mathbf Z = \{\lambda \in \mathbf R:$
il existe
$\lambda_i$
tel que
$\lambda \equiv \lambda_i \bmod \mathbf Z\}.$
Les autres notations sont celles de l'\'enonc\'e $\bf 2$.
\label{thm52}
\end{thm}

{\bf Premi\`ere partie de la d\'emonstration: Preuve des \'enonc\'es 1 et 2 pour $P_{k,l},$ 
$ G_{2m+6},$ $ G_{2m+3},$ $U_{2m+1},$ $ FW_{13},$ $FW_{19}, K_{13}.$}

On se rappelle que l'argument de la d\'emonstration du Th\'eor\`eme ~\ref{thm1}
~\ref{thm43} s'appuyait sur $\bf (a)$ le fait que 
les \'el\'ements des matrices $P^{(1)},P^{(2)}$
sont monomiaux, $\bf (b)$ apr\`es un changement de base de
$\tilde \Phi$
et $F$, on peut supposer que
$P^{(1)}(s)$ est 
une matrice diagonale
et que $P^{(2)} =$
(matrice diagonale)
$\times \Sigma$
avec
$\Sigma \in SL(\mu, {\mathbf Z})$ telle que $ \Sigma^{\mu} = id_{\mu}$.

Les \'enonc\'es {\bf 1} et {\bf 2} se d\'eduisent de $\bf (a)$ et
$\bf (b).$ 

{\it Preuve de \bf (a)} 
D'abord on se souvient de la s\'erie de Poincar\'e $P_V(t), (1.2.3)$
dans le  cas $n=1, k=2, m=3,$
$$    P_{V}(t)    =          
(1-t^{p_1})(1-t^{p_2})\frac{(t^{w_1}+ t^{w_2}+ t^{w_3}-t^{p_1}- 
t^{p_2}-1)} {(1-t^{w_1})(1-t^{w_2})(1-t^{w_3})} +1. \leqno(5.2.3) $$ 
Cela et la relation (1.5.1) donnent: 
$$ w (P^{(1)}_{ij}(s)) = p_2 + w(\tilde{\omega_i}) - w(\phi_j)- |w| \leq 2 
(|p|-|w|)- w(\phi_j). \leqno(5.2.4)$$
D'autre part, pour que un \'el\'ement de $ P^{(1)}_{ij}(s)$ consiste de 
plusieurs mon\^omes, 
il faut qu'il existe des paires de vecteurs entiers $(\alpha, \beta)\not = 
(\alpha', \beta') \in {\bf Z}^2_{\geq 0}$ telles que:
$$w (P^{(1)}_{ij}(s)) = \alpha p_1 + \beta p_2 = {\alpha}' p_1 + {\beta}' p_2.$$
Pour les singularit\'es de 5.1, 5.2 ci-dessus, il est impossible de 
trouver de telles paires de vecteurs 
entiers positifs  d'apr\`es (5.2.4). Quant \`a la matrice
$P^{(2)}$ le m\^eme argument fonctionne.

{\it Preuve de \bf (b).}

Le lemme ~\ref{lem8} a d\'ej\`a d\'emontr\'e 
la propri\'et\'e $\bf (b)$
dans les cas 
non-r\'esonants. 

Afin de le prouver dans les cas unimodaux qui figurent dans l'\'enonnc\'e 
{\bf 2}, on 
reproduit 
l'argument du lemme ~\ref{lem8}. Plus pr\'ecisement, il suffit de 
d\'emontrer la validit\'e de la  d\'ecomposition comme (5.1.8) 
 \`a laquelle en fait
un seul terme prend part. 
On va d\'emontrer que si on a une d\'ecomposition sous la forme suivante,
$$ \tilde {\omega_j} \wedge df_2= (P^{(1)}_{jj_2}\phi_{j_2} 
+ P^{(1)}_{jj_2'}\phi_{j_2'})dx, \leqno(5.2.5)$$
$$ \tilde {\omega_j} \wedge df_1= (P^{(2)}_{jj_1}\phi_{j_1} 
+ P^{(1)}_{jj_1'}\phi_{j_1'})dx,$$
alors $P^{(1)}_{jj_2} =0$ ou $P^{(1)}_{jj_2'}=0$ 
($P^{(2)}_{jj_1} =0$ ou $P^{(2)}_{jj_1'}=0$).
Il est \'evident qu'il y aurait  au plus deux termes \`a 
droite de (5.2.5), car   $ \tilde {\omega_j}$ consiste au plus en 
deux termes.

La liste ci-dessus donne la d\'emonstration de l'\'enonc\'e $\bf (b)$
pour les cas
$P_{k,l}, G_{2m+3}, G_{2m+6}, U_{2m+1}.$

Quant aux cas
$FW_{13}, FW_{19}, K_{13}$
il suffit de comparer les poids quasihomog\`enes des formes
$\tilde \omega_j \in F$
et ceux de
$\phi_i(x) dx, \phi_i(x) \in \tilde\Phi$.
(Il n'est pas n\'ec\'essaire de calculer toutes les formes de $F$):

\begin{flushleft}
$FW_{3r+1} (r = 4,6) :$
\end{flushleft}

\[
\left\{
\begin{array}{ccccc}
f_1 & = & x_1x_3 + x^r_2 \\ 
f_2 & = & x_1x_2 + x^3_3 
\end{array}
\right.
\]
$$\begin{array}{c}
w(\tilde \omega_j) = \{ r+1+4(i+1), 2(r+1)+4(i+1), 3(r+1)+4i (0 \le i \le r-1), 7r+3\} \\
w(\phi_j(x)dx) = \{ 4(r+1)+4i, 5(r+1)+4i, 6(r+1)+4i (0 \le i \le r-1), 7r+3\}
\end{array}$$
\begin{flushleft}
$K_{13}:$
\end{flushleft}

\[
\left\{
\begin{array}{c}
f_1 = x^2_1 + x^3_2 x_3  \\
f_2 = x_1x_2 + x^2_3
\end{array}
\right.
\]
$$\begin{array}{c}
w(\tilde \omega_j)  =  \{8,10,11,12,13,14,15,16,17,18,19,20,23\} \\
w(\phi_j(x)dx)  =  \{15,18,20,21,22,23,24,25,26,27,28,30,33\}
\end{array}$$
Pour que la situation (5.2.5) avec deux termes non-nuls 
se produise, il est n\'ec\'essaire qu'il existe
$\phi_{j_1}, \phi_{j_1'},\phi_{j_2},\phi_{j_2'} \in \tilde \Phi;$ 
$\alpha, \beta,$
$\gamma, \delta,$ $\alpha^{\prime}, \beta^{\prime},$ 
$\gamma^{\prime}, \delta^{\prime}$ $\in 
{\mathbf Z_{\geq 0}},$ $j_1 \ne j^{\prime}_1, j_2 \ne j^{\prime}_2$
tels que
$$\begin{array}{c}
w(\tilde \omega_j) + p_2  =  w(\phi_{j_2}) +\alpha p_1 + \beta p_2+|w| 
= w(\phi_{j^{\prime}_2}) + \alpha^{\prime}p_1 + \beta^{\prime}p_2 +|w|\\
w(\tilde \omega_j) + p_1  =  w(\phi_{j_1}) + \gamma p_1 + \delta p_2 +|w|
= w(\phi_{j^{\prime}_1}) + \gamma^{\prime}p_1 + \delta^{\prime}p_2+|w|.
\end{array} \leqno{(5.2.6)}$$
Afin de voir $\bf(b)$, il suffit de le v\'erifier pour
$\tilde \omega_j \in F$
qui satisfait (5.2.6). Il n'existe pas de telle forme
$\tilde \omega_j \in F$
dans les cas
$FW_{13}, FW_{19}$.
Dans le  cas
$K_{13}$,
on trouve les formes satisfaisant (5.2.6) comme suit:
$$\{\tilde \omega_{15}= 3x_1dx_2 dx_3 + 2x_2dx_3dx_1, \tilde \omega_{18} = x_2 \tilde \omega_{15}, \tilde \omega_{20} = x_3 \tilde \omega_{15}, \tilde \omega_{23} = x_2 x_3 \tilde \omega_{15}\}$$
o\`u l'indice  de chaque forme indique 
son poids quasihomog\`ene. On voit facilement que la situation (5.2.5)
avec deux termes non-nuls 
n'arrive  pour aucune des formes ci-dessus. On a donc d\'emontr\'e que 
dans les 
cas
de $FW_{13}, FW_{19}, K_{13}$,
on trouve au plus des entiers uniques
$j_1, j_2 \in [1,\mu]$
pour chaque
$j \in [1,\mu]$
tels que
$$\begin{array}{c}
\tilde \omega_{j_2} \wedge df_1  =  P^{(2)}_{j_2j} (f) \phi_j(x) dx \\
\tilde \omega_{j_1} \wedge df_2  =  P^{(1)}_{j_1j} (f) \phi_j(x)dx 
\end{array}$$
avec
$P^{(2)}_{j_2j} \ne 0, P^{(1)}_{j_1j} \ne 0$. L'\'enonc\'e  {\bf (b)} s'en 
d\'eduit imm\'ediatement.

{\bf Seconde partie: Preuve des \'enonc\'es 1 et 3  pour les cas
$FZ_{6m+6},$ $FZ_{6m+8},$ $ HD_{13},$ $ HD_{14},$ $ K_{14}$.}

Dans ces cas, on trouve des formes
$\phi_j(x)dx, \phi_j \in \tilde \Phi$
et
$1 \le i \ne i^{\prime} \le \mu$
tels que

$$\begin{array}{c}
\tilde \omega_i \wedge df_l  =  P^{(l)}_{ij} \phi_j(x)dx \\
\tilde \omega_{i^{\prime}} \wedge df_l = P^{(l)}_{i^{\prime}j} \phi_j(x)dx
\end{array}$$
avec $P^{(l)}_{ij}\ne 0, P^{(l)}_{i^{\prime}j}\ne 0$ et
$l \in \{1,2\}$.
Cela implique que $\bf (b)$ 
n'est pas valable pour ces singularit\'es. Pourtant on trouve parmi 
elles des sous-ensembles
$\tilde \Phi_0 := \{\phi_{j_1}, \ldots, \phi_{j_\nu}\} \subset \tilde \Phi$
et
$F_0 := \{ \tilde \omega_{i_1}, \ldots, \tilde \omega_{i_\nu}\} \subset F, 
2 \le \nu \le \mu,$ tels que les matrices
$P^{(1)}_0(s), P^{(2)}_0(s) \in GL(\nu, \mathbf C) \otimes \mathbf C[s]$
d\'efinies comme dans (1.5.1) pour les formes de
$F_0$
et
$\phi_{j_1}dx, \ldots, \phi_{j_\nu} dx$
soient toutes les deux des produits d'une matrice 
diagonale et d'une matrice de permutation. Bien entendu, 
pour la transform\'ee de Mellin
$M_{{j_k}}(z)$
d\'efinie par le mon\^ome
$\phi_{j_k} \in \tilde \Phi_0$
un \'enonc\'e parallel \`a celui du Th\'eor\`eme 3.1 est valable.

Pour \'etablir l' \'enonc\'e $\bf 1$, on reproduit un argument similaire \`a 
celui utilis\'e pour voir $\bf (b)$ dans les cas
$FW_{13}, FW_{19}, K_{13}$. Les formes satisfaisant (5.2.6) sont \'enum\'er\'ees
ci-dessous. On \'etablit aussi une liste des poids quasihomog\`enes de
$\tilde \omega_w \in F$,
l'espace
$\tilde \Phi$
et des poids quasihomog\`enes
$\phi_j(x)dx, \phi_j(x) \in \tilde \Phi$.
\[
$$\begin{eqnarray*}
HD_{13} :f_1 &= &x^2_1 + x^2_2x_3\\
 f_2& = &x_1x_2 + x^3_3\\ 
w(\tilde \omega_j) & = & \{9,11,12,13,14,15,16,17,18,19,20,22,24\} \\
\tilde \Phi & = & \{1,x_3,x_2,x_1,x^2_3, x_2x_3, x^2_2, x_1x_3, x_1x_2, x_1x^2_3, 
x_1x_2x_3,x_1x_2x^2_3\} \\
w(\phi_jdx) & = & \{16,20,21,23,24,25,26,27,28,29,31,32,36\} \\
\tilde \omega_{12} & = & dx_1dx_2 \\
\tilde \omega_{14} & = & x_2dx_2dx_3 \\
\tilde \omega_{16} & = & x_1 dx_2dx_3 + x_2dx_3dx_1 \\
\tilde \omega_{18} & = & x_2x_3 dx_2 dx_3 + x_1 dx_3 dx_1 \\
\tilde \omega_{20} & = & 3x_1x_3 dx_2dx_3 + x^2_3 dx_1 dx_2 \\
\tilde \omega_{24} & = & x^2_3 \tilde \omega_{16} 
\end{eqnarray*}$$
\]

\[
$$\begin{eqnarray*}
HD_{14} : f_1 & = & x^2_1+x^3_2 \\
f_2 & = & x_1x_2 + x^3_3 \\
\tilde \omega_{17} & = & x_2dx_2dx_3 \\
\tilde \omega_{20} & = & 3x_1dx_2dx_3 + 2x_2 dx_3 dx_1 \\
\tilde \omega_{26} & = & x_2 \tilde \omega_{20} \\
\tilde \omega_{23} & = & x_1dx_3dx_1 + x^2_2 dx_2 dx_3 \\
w(\tilde \omega_j) & = & \{11,14,15,16,17,19,20,21,22,23,25,26,28,31\} \\
\tilde \Phi & = & \{1,x_3,x_2,x_1,x^2_3, x_2x_3, x^2_2, x_1x_3, x_1x_2, x_2x^2_3, x^2_2x_3, x^4_3, x_1x^2_2, x_2x^4_3\} \\
w(\phi_jdx) & = & \{20,25,26,29,30,31,32,34,35,36,37,40,41,46\}
\end{eqnarray*}$$
\]

\[
$$\begin{eqnarray*}
K_{14}: f_1 & = & x^2_1 + x^5_2 \\
f_2 & = & x_1x_2 + x^2_3 \\
w(\tilde \omega_j) & = & \{11,14,15,17,18,19,21,22,23,25,26,27,29,33\} \\
\tilde \Phi & = & \{1,x_2,x_3, x^2_2, x_1, x_2x_3, x^3_2, x_1x_2, x^2_2x_3, x^4_2, x_1x^2_2, x^3_2x_3, x_1x^3_2, x_1x^4_2\} \\
w(\phi_jdx) & = & \{21,25,28,29,31,32,33,35,36,37,39,40,43,47\} \\
\tilde \omega_{15} & = & x_2dx_2dx_3 \\
\tilde \omega_{19} & = & x_2\tilde \omega_{15} \\
\tilde \omega_{22} & = & x^2_2 dx_1dx_2 \\
\tilde \omega_{21} & = & 5x_1dx_2dx_3 + 2x_2dx_3dx_1 \\
\tilde \omega_{25} & = & x_2 \tilde \omega_{21} \\
\tilde \omega_{23} & = & x^3_2dx_2 dx_3 \\
\tilde \omega_{27} & = & x_1dx_3dx_1 + x^4_2dx_2 dx_3 \\
\tilde \omega_{29} & = & x^2_2 \tilde \omega_{21} \\
\tilde \omega_{33} & = & x^3_2 \tilde \omega_{21} 
\end{eqnarray*}$$
\]
Ici l'indice de chaque forme indique son poids quasihomog\`ene. Le calcul direct montre que la situation (5.2.5) avec deux termes non-nuls
ne peut se produire  pour aucune des formes 
ci-dessus. Cela d\'emontre l'\'enonc\'e {\bf 1} pour ces singularit\'es.

Quant aux mon\^omes de
$\tilde \Phi \setminus \tilde \Phi_0$,
apr\`es un calcul simple, on observe la situation suivante. Pour les mon\^omes 
$\phi_{j_{\nu+l(r)}} \in \tilde \Phi \setminus \tilde \Phi_0, 1 \le r 
\le \mu-\nu$
il existe au moins un mon\^ome
$\phi_{j_{\nu+l(0)}} \in \tilde \Phi_0$
et des mon\^omes
$\phi_{j_{\nu+l(1)}}, \phi_{j_{\nu+l(2)}}, \ldots, \phi_{j_{\nu+l(r-1)}} \in 
\tilde \Phi \setminus \tilde \Phi_0$
tels que
$$\begin{array}{c}
(z_1+\tilde{\eta}_{j_{\nu+l(q)}} + \lambda_{j_{\nu+l(q)}}+1) M_{j_{\nu+l(q)}} (z_1 + \tilde{\eta}_{j_{\nu+l(q)}}, z_2) \\
= \tilde v_{j_{\nu+l(q)}}z_1 M_{j_{\nu+l(q-1)}} (z_1-1, z_2 + \tilde{\delta}_{j_{\nu+l(q)}} + 1), q = 1,\ldots, r
\end{array} \leqno{(5.2.7)}$$
pour certains
$\tilde{\eta}_{j_{\nu+l(q)}}, \tilde{\delta}_{j_{\nu+l(q)}} \in 
{\mathbf Z}_{\ge 0}$
et
$\tilde v_{j_{\nu+l(q)}} \in \mathbf Q$.
En contraste avec le cas des mon\^omes de
$\tilde\Phi_0$, la relation de r\'ecurrence (5.2.7) ne se ferme pas pour
$\phi_{j_{\nu+l(q)}} \in \tilde \Phi \setminus \tilde \Phi_0$. Voir (4.1.1).
C'est-\`a-dire qu'il n'existe pas d'indices
$0 \le r^{\prime\prime} < r^{\prime} \le r$
tels que une
$(r+1)$
\`eme relation de r\'ecurrence comme (5.2.7) soit v\'erifi\'ee pour
$M_{j_{\nu+l(r^{\prime})}}$
et
$M_{j_{\nu+l(r^{\prime\prime})}}, r^{\prime\prime} \ne r^{\prime}-1$.
On note des mon\^omes de
$\tilde \Phi \setminus \tilde \Phi_0$
ci-dessous:

\[
$$\begin{eqnarray*}
FZ_{6m+6}, FZ_{6m+8} & : & \tilde \Phi \setminus \tilde \Phi_0 = \{x^2_3\} \\
HD_{13} & : & \tilde \Phi \setminus \tilde \Phi_0 = \{x_1x_2x_3\} \\
HD_{14} & : & \tilde \Phi \setminus \tilde \Phi_0 = \{x_2x_3, x_1x_3, x_1x_2, x^2_2x_3, x^4_3, x_1x^2_2, x_2x^4_3\} \\
K_{14} & : & \tilde \Phi \setminus \tilde \Phi_0 = \{x_2, x^2_2x_3, x^3_2x_3\} 
\end{eqnarray*}$$
\]
En reproduisant l'argument du Th\'eor\`eme 4.3 il est facile de voir que
$$\begin{array}{c}
M_{j_{\nu+l(0)},\gamma_q} (z) = 
V^{\frac{z_1}{d_1(j_{\nu +l(0)})}}_{j_{\nu+l(0)}} 
\prod_{j\in I^{(j_{\nu+l(0)})}_1} 
\Gamma(-\frac 1 {d_1(j_{\nu+l(0)})} L(\ast, j_{\nu+l(0)}, j; z_1)) \\
\times \prod_{i \in I^{(j_{\nu+l(0)})}_2} 
\Gamma (-\frac {\tilde L(\ast, j_{\nu+l(0)},i;z_1)}{d_1(j_{\nu+l(0)})}) 
\prod_{t\in I^{(j_{\nu+l(0)})}_1 \cup I_2^{(j_{\nu+l(0)})}} 
\Gamma(-\frac{(z_1+t)}{d_1(j_{\nu+l(0)})})^{-1} \\
(p_2(z_2-\zeta_q) + p_1z_1)^{-1}g(z_1)
\end{array} \leqno{(5.2.8)}$$
avec les notations de l'\'enonc\'e $\bf 2$. 
L'expression (5.2.8) et la relation de r\'ecurrence pour
$\phi_{j_{\nu+l(q)}}, q = 0,1, \ldots, r$
entra\^{\i}nent le r\'esultat.
\hfill C.Q.F.D.

\begin{remark}
{\em
D'apr\`es la m\'ethode pr\'esent\'ee ci-dessus au $\S$ 4, on n'est pas susceptible de r\'esoudre le syst\`eme de Gauss-Manin associ\'e aux SIIC unimodales quasihomog\`enes de courbe espace qui ne figurent pas dans le Th\'eor\`eme 5.2.

Pour ces singularit\'es i.e.
$FT_{4,4}, FW_{14}, FW_{1,0}, FW_{18}, FZ_{6m+7}, 
FZ_{6m-1,0}, J_{6m+8}, J_{m+1,0}, K_{1,0}$
de la liste de Wall, il existe au moins une forme
$\tilde \omega_j \in F$
qui suscite la situation (5.2.5) avec deux termes non-nuls. 
Par example, pour la singularit\'e
$FT_{4,4}$:

 \[ 
\left\{ 
\begin{array}{ccc} 
f_1 &=& x_1x_2 + x^3_3 \\
f_2 &= & x_1x_3 + x^3_2 + \lambda x_2x^2_3, \lambda \ne 0.
\end{array} 
\right.
\] 

On trouve trois formes
$\tilde \omega_6, \tilde \omega_7, \tilde \omega_8$
du m\^eme poids quasihomog\`ene i.e.

\[
$$\begin{eqnarray*}
w(\tilde \omega_6) & = & w(\tilde \omega_7) = w(\tilde \omega_8) = 4, \qquad\mbox{car}\qquad w_1 = 2, w_2 = w_3 =1 \\
\tilde \omega_6 & = & 2x_2 dx_1 dx_2 + (-4\lambda x^2_2 + (2\lambda^2 +6)x^2_3) dx_2dx_3 \\
df_1 \wedge \tilde \omega_6 & = & \lambda (-4f_1 + x_1x_3)dx \\
df_2 \wedge \tilde \omega_6 & = & (2f_1 -8x^3_3)dx \\
\tilde \omega_7 & = & x_1dx_2 dx_3 \\
df_1 \wedge \tilde \omega_7 & = & (f_1 - x^3_3) dx \\
df_2 \wedge \tilde \omega_7 & = & x_1x_3 dx \\
\tilde \omega_8 & = & 2x_2dx_1dx_2 - \lambda x_3dx_3dx_1 - (\lambda x^2_2 + (\lambda^2+6)x^2_3) dx_2dx_3 \\
df_1 \wedge \tilde \omega_8 & = & -\lambda f_2dx \\
df_2 \wedge \tilde \omega_8 & = &( 2f_1-(8+2\lambda^2) x^3_3)dx \\
\tilde \Phi & = & \{1,x_2, x_3, x_1, x^2_2, x^2_3, x_2x_3, x_1x_3, x^3_3, x_1x^2_2\}.
\end{eqnarray*}$$
\]
On peut voir qu'il est impossible d'\'eviter la situation (5.2.5)
avec deux termes non-nuls en 
choisissant une autre base de $F$ et
$\Phi$.

Les espaces vectoriels
$F,\Phi$
et les matrices
$P^{(1)}(s), P^{(2)}(s)$
sont disponibles sur demande pour toutes les SIIC unimodales quasihomog\`enes de courbe espace.}
\label{remark52}
\end{remark}
\begin{remark}

{\em Tout r\'ecemment, l'auteur a reussi \`a obtenir le r\'esultat suivant.

Les spectres du syst\`eme (1.10.1), d\'efinis comme dans la D\'efinition
~\ref{dfn411} consistent des donn\'ees du type;
 \[ 
\left\{ 
\begin{array}{ccc} 
p_1z_1 + p_2 z_2 + w(\psi_j)+ p_2-p_1& \leq & 0\\ 
z_1 \leq -1& \rm{ou} & z_2 \leq -1. \\ 
\end{array} 
\right. 
\] 
 De cet \'enonc\'e, il suit facilement que les spectres du syst\`eme de 
Gauss-Manin 
associ\'e \`a une ICIS quasihomog\`ene de courbe espace sont sym\'etriques
par rapport \`a une droite. On peut ainsi en d\'eduire que le Th\'eor\`eme
~\ref{thm4} ci-dessous est valable pour ces singularit\'es, si on pose
$ \lambda_i = \frac{1}{p_1} w(\psi_i), \; 1 \leq i \leq \mu.$

La d\'emonstration s'appuie sur une m\'ethode appel\'ee ``Cayley trick''
qui calcule la cohomologie relative de de Rham d'une SIIC (0.1) au moyen de la cohomologie d'une hypersurface d\'efinie comme suit:
$$ X_f =\{(x_1, \cdots,x_{n+k} , y_1, \cdots, y_k) \in {\bf C}^{n+2k} ; 
y_1(f_1(x)-s_1)+ \ldots y_k(f_k(x) -s_k) =0 
\}.$$ Le detail doit para\^{i}tre dans un travail suivant.} 
\label{remark53}
\end{remark}
\vspace{2pc} 
{ 
\center{\section{ 
Nombre de Hodge des fibres de Milnor }} 
} 
On revient \`a la situation des chapitres $\S 0,$ $\S 1.$ 
Soit $X_s$ une fibre de Milnor de SIIC courbe espace pour 
$(s_1, s_2) \in S$ hors sa valeur critique.  On 
regarde une vari\'et\'e projective ${\bar X}_s \subset {\bf P} 
(w_1, w_2,w_3,1)$
associ\'ee \`a $ X_s :$
$${\bar X}_s :=\{(x_1,x_2,x_3) \in {\bf P} (w_1, w_2,w_3,1) ; f_1(x_1,x_2,x_3) = 
s_1 u^{p_1}, f_2(x_1,x_2,x_3) = s_2 u^{p_2}\}.$$
C'est une cl\^oture de ${\bar X}_s$ dans l'espace projectif 
${\bf P} (w_1, w_2,w_3,1)$
avec des poids comme en (1.1.1). On prend le poids de la variable $u$ 
\'egal \`a $1.$  On peut consulter \cite{Dolg}, \cite{Dim} pour la 
projectivisation 
d'une SIIC.

On note $p_g(\bar{X}_s)$ le genre g\'eom\'etrique de la courbe 
${\bar X}_s$ qui doit \'egaler 
le rang du groupe de cohomologie $H^1({\bar X}_s, {\cal O}_{{\bar X}_s} ).$ Le 
genre g\'eom\'etrique de la fibre de Milnor
est un des invariants importants de la singularit\'e 
$X_0.$  
On a une expression assez simple de
cet invariant au moyen des spectres du syst\`eme de Gauss-Manin.

Nous nous  rappelons  ici que les nombres de Hodge $h^{pq}(X_s) = Gr_F^{p}Gr^W_{p+q} H^n(X_s)$
(et sa s\'erie de Poincar\'e) eux-m\^emes ont \'et\'e  calcul\'es 
par F.Hirzebruch pour le cas d'une SIIC homog\`ene 
 \cite{Hir}, \S 22 et par H.Hamm pour le cas d'une SIIC 
quasihomog\`ene de dimension positive quelconque \cite{Ham}.

 \begin{thm}
Pour les singularit\'es trait\'ees dans les Th\'eor\`emes 4.3 et 5.2, 
on a les formules  suivantes
$$ h^{01}(X_s) = h^{10}(X_s) = p_g({\bar X}_s) = \sharp \{ i ; \lambda_i < 0 \} \leqno(6.1)$$
$$ h^{11}(X_s) = p_g(X_0) = \sharp \{ i ; \lambda_i = 0 \} \leqno(6.2)$$
o\`u $\lambda_i, 1 \leq i \leq \mu $ sont les spectres du syst\`eme de Gauss-Manin (3.1.1),
(5.2.1).

\label{thm4}
\end{thm}

{\bf D\'emonstration}
{\bf i)} D\'emonstration de (6.1).
Tout d'abord, rappelons que les  ${\bf C}-$modules  $A_{X_0}, A_{{\bar X}_s}$
des polyn\^omes sur  ${X_0}$  et ${{\bar X}_s}$ sont:
$$A_{X_0} := {\bf C}[x_1,x_2,x_3]/(f_1, f_2) , 
A_{{\bar X}_s} = {\bf C}[x_1,x_2,x_3,u]/
(f_1 -s_1 u^{p_1}, f_2-s_2 u^{p_2}). $$

Ces deux modules $A_{X_0}, A_{{\bar X}_s}$ sont munis d'une 
filtration naturelle
compatible avec le poids $w(g)$ de $g \in A_{X_0}, A_{{\bar X}_s}.$ Ici
$$ ( E + u \frac{\partial}{\partial u}) g(x_1,x_2,x_3,u) = w(g) g(x_1,x_2,x_3,u)$$
avec le champ de Euler introduit par (1.1.2). On peut regarder les s\'eries
de Poincar\'e $P_{A_{X_0}}(t), P_{A_{{\bar X}_s}}(t)$ d\'efinies par une 
filtration 
de poids :
$$P_{A_{X_0}}(t)= \sum_{ a_d = \sharp\{ g(x_1,x_2,x_3) \in A_{X_0}; w(g)=d \}} 
a_d t^d,$$
$$ P_{A_{{\bar X}_s}}(t)= \sum_{ a_d = \sharp 
\{ h(x_1,x_2,x_3,u) \in A_{{\bar X}_s}; w(h)=d \}} a_d t^d.$$
Selon Dolgachev \cite{Dolg} 3.4.4, elles satisfont la relation suivante:
$$P_{A_{X_0}}(t) = (1 -t)P_{A_{{\bar X}_s}}(t)$$
o\`u
$P_{A_{X_0}}(t) =    \frac{(1-t^{p_1})(1-t^{p_2})} {(1-t^{w_1})(1-t^{w_2})(1-t^{w_3})}.$
Par contre on a 
$$    P_{\tilde \Phi}(t)    =    t^{|p|- |w|} + (1-t^{|p|- |w|})P_{A_{X_0}}(t)  \leqno(6.3)$$
d'apr\`es la formule d'Aleksandrov (3.4.6). Ici $|w|= w_1 + w_2 + w_3$
et $|p| = p_1 + p_2.$  En somme
$$P_{A_{{\bar X}_s}}(t) = 
\frac{  P_{\tilde \Phi}(t)-t^{|p|- |w|}}{(1-t)(1-t^{|p| - |w|})} 
= (P_{\tilde \Phi}(t)-t^{|p|- |w|})(\sum_{j \geq 0} t^j)
(\sum_{k \geq 0}t^{k(|p|- |w|)}). \leqno(6.4)$$
Le th\'eor\`eme de Dolgachev cit\'e plus haut implique que le  
coefficient de $t^{|p|- |w|-1}$ de $P_{A_{{\bar X}_s}}(t)$ 
donne le genre g\'eom\'etrique $p_g({\bar X}_s).$
Si on note le d\'eveloppement de $P_{\tilde \Phi}(t):$ 
$$    P_{\tilde \Phi}(t) = \sum_{j \geq 0} \pi_j t^j, $$
on d\'eduit de la relation $(6.4), $
$$p_g({\bar X}_s) = 2	\pi_0 + \pi_1 + \cdots \pi_{|p|- |w|-1}-1,$$
$$= \sharp\{ \mbox{polynomes de}\; \tilde \Phi \;\mbox{de poids}\leq |p|- |w|-1  \}$$
$$= \sharp\{ \mbox{polynomes de} \; \tilde \Phi \; \mbox{de poids} < |p|- |w|  \}.$$
L'unicit\'e de l'\'element de poids $0$ dans $\tilde \Phi $
entra\^{\i}ne la premi\`ere \'egalit\'e. 
La derni\`ere egalit\'e donne (6.1) en utilisant les Th\'eor\`emes ~\ref{thm1} et ~\ref{thm52}.

{\bf ii)} D\'emonstration de (6.2). Quant \`a la seconde \'egalit\'e, on 
reproduit l'argument 
bien utilis\'e depuis Steenbrink \cite{St1}. Voir aussi  
\cite {Dim}, \cite{Al1}, (5.4). 
Soit $Y = {\bar X}_s \setminus X_s.$
A partir de la suite exacte
 $$ \cdots \rightarrow H^{n-2}(Y)(-1) \rightarrow H^n({\bar X}_s) \rightarrow H^n({X}_s)
\rightarrow H^{n-1}(Y)(-1)\rightarrow H^{n+1}({\bar X}_s)\rightarrow \cdots, $$
on obtient une suite exacte courte 
$$0 \rightarrow P^n({\bar X}_s) \rightarrow H^n({X}_s)
\rightarrow P^{n-1}(Y)(-1)\rightarrow,\leqno(6.5) $$
o\`u 
 $$P^n({\bar X}_s) = coker( H^{n-2}(Y)(-1)\rightarrow H^n({\bar X}_s))$$
$$P^{n-1}(Y)(-1) = ker (H^{n-1}(Y)(-1) \rightarrow H^{n+1}({\bar X}_s))$$
qui s'appellent la partie primitive de la cohomologie correspondente.
De (6.5), on obtient 
$$h^{p,n+1-p}(X_s) = Gr_F^{p}Gr^W_{n+1} H^n(X_s)= Gr_F^p P^{n-1}(Y)(-1) = h_0^{p-1,n-p}(Y).$$
Le dernier est calcul\'e par le Th\'eor\`eme 4.4, 3) de \cite{Al1} qui dit que 
$$ h_0^{0,0}(Y) = \mbox{ coefficient de}\; t^{|p|- |w|}\; \rm{de} \; P_{A_{X_0}}(t),$$
qui est \'egal \`a son tour au
coefficient de $t^{|p|- |w|}$  de $P_{\tilde \Phi}(t),$
par (6.3). 
La d\'emonstration s'ach\`eve si on se souvient de la d\'efinition des 
spectres $\lambda_i$ des Th\'eor\`emes ~\ref{thm1}  et ~\ref{thm52}.
{\bf C.Q.F.D.}


\vspace{\fill} 

%

\noindent

\begin{flushleft} 
\begin{minipage}[t]{6.2cm} 
  \begin{center} 
{\footnotesize Institute of Control Sciences\\ 
Russian Academy of Sciences,\\ 
Profsojunaja ul. 65,\\ 
GSP-7, Moscow, 117806,\\ 
Russia\\ 
{\it E-mails}:  tanabe@ipu.rssi.ru\\
\vspace{1pc}
Max Planck Institut f\"ur Mathematik\\ 
Vivatsgasse 7,Bonn,D- 53111, Germany\\
{\it E-mail:} tanabe@mpim-bonn.mpg.de}
\end{center} 
\end{minipage} 
\end{flushleft}  
 
\end{document}